\documentclass[12pt]{article}

\title{The local Gromov-Witten theory of curves}
\author{Jim Bryan and Rahul Pandharipande}

\usepackage{diagrams}

\usepackage{texdraw}
\input{txdtools}

\def\fatness{1.2}
\everytexdraw{%
	\drawdim pt \linewd 0.5 \textref h:C v:C
  \arrowheadsize l:4 w:3  
  \arrowheadtype t:V
}

\setbox7=\hbox{
\raisebox{1pt}{\begin{texdraw}
	\drawdim pt \linewd 0.3 \setunitscale 1.2
	\move (-3 0)
	\move (0 -1.2) 
	\rlvec (6 0) 
	\rlvec (0 -1.2)
	\rlvec (2.4 2.4) 
	\rlvec (-2.4 2.4)
	\rlvec (0 -1.2) 
	\rlvec (-6 0)
	\rlvec (0 -2.4)
	\move (11.5 0)
\end{texdraw}}}

\newcommand{\idFA}{%
	\bsegment
		\realmult {\fatness} {2.5} \xrad
		\realmult {\fatness} {6.25} \yrad
		\move (0 0) \lellip rx:{\xrad} ry:{\yrad}
		\move (40 0) \lellip rx:{\xrad} ry:{\yrad}
		\realadd { 0} { \yrad} \inB
		\realadd { 0} {-\yrad} \inA
		\move (0 {\inB}) \rlvec (40 0) 
		\move (0 {\inA}) \rlvec (40 0) 
		\savepos (40 0)(*ex *ey)
	\esegment
	\move (*ex *ey)
}
\newcommand{\counitFA}{%
	\bsegment
		\realmult {\fatness} {2.5} \xrad
		\realmult {\fatness} {6.25} \yrad
		\move (0 0)
		\move (20 0)
		\move (0 0) \lellip rx:{\xrad} ry:{\yrad}
		\realadd { 0} { \yrad} \inB
		\realadd { 0} {-\yrad} \inA
		\move (0 {\inB}) \rlvec (8 0) 
		\move (0 {\inA}) \rlvec (8 0) 
	  \move (8 {\inA}) \clvec (20 {\inA})(20 {\inB})(8 {\inB}) 

		\savepos (0 0)(*ex *ey)
	\esegment
	\move (*ex *ey)
}
\newcommand{\unitFA}{%
	\bsegment
		\realmult {\fatness} {2.5} \xrad
		\realmult {\fatness} {6.25} \yrad
		\move (0 0)
		\move (-20 0)
		\move (0 0) \lellip rx:{\xrad} ry:{\yrad}
		\realadd { 0} { \yrad} \inB
		\realadd { 0} {-\yrad} \inA
		\move (0 {\inB}) \rlvec (-8 0) 
		\move (0 {\inA}) \rlvec (-8 0) 
	  \move (-8 {\inA}) \clvec (-20 {\inA})(-20 {\inB})(-8 {\inB}) 

		\savepos (0 0)(*ex *ey)
	\esegment
	\move (*ex *ey)
}
\newcommand{\twistFA}{%
	\bsegment
		\realmult {\fatness} {2.5} \xrad
		\realmult {\fatness} {6.25} \yrad
		\move (0 0) \lellip rx:{\xrad} ry:{\yrad}
		\move (0 30) \lellip rx:{\xrad} ry:{\yrad}
		\move (40 0) \lellip rx:{\xrad} ry:{\yrad}
		\move (40 30) \lellip rx:{\xrad} ry:{\yrad}
		\realadd {30} { \yrad} \inD
		\realadd {30} {-\yrad} \inC
		\realadd { 0} { \yrad} \inB
		\realadd { 0} {-\yrad} \inA
		\move (0 {\inD}) \clvec (28 {\inD})(20 {\inB})(40 {\inB}) 
		\move (0 {\inC}) \clvec (20 {\inC})(12 {\inA})(40 {\inA}) 

		\move (0 {\inB}) \clvec (20 {\inB})(12 {\inD})(40 {\inD}) 
		\move (0 {\inA}) \clvec (28 {\inA})(20 {\inC})(40 {\inC}) 
		\savepos (40 0)(*ex *ey)
	\esegment
	\move (*ex *ey)
}

\newcommand{\multFA}{%
	\bsegment
		\realmult {\fatness} {2.5} \xrad
		\realmult {\fatness} {6.25} \yrad
		\move (0 0) \lellip rx:{\xrad} ry:{\yrad}
		\move (0 30) \lellip rx:{\xrad} ry:{\yrad}
		\move (40 15) \lellip rx:{\xrad} ry:{\yrad}
		\realadd {15} { \yrad} \outB
		\realadd {15} {-\yrad} \outA
		\realadd {30} { \yrad} \inD
		\realadd {30} {-\yrad} \inC
		\realadd { 0} { \yrad} \inB
		\realadd { 0} {-\yrad} \inA
		\move (0 {\inD}) \clvec (25 {\inD})(20 {\outB})(40 {\outB}) 
		\move (0 {\inB}) \clvec (18 {\inB})(18 {\inC})(0 {\inC}) 
		\move (0 {\inA}) \clvec (25 {\inA})(20 {\outA})(40 {\outA}) 
		\savepos (40 15)(*ex *ey)
	\esegment
	\move (*ex *ey)
}

\newcommand{\comultFA}{%
	\bsegment
		\realmult {\fatness} {2.5} \xrad
		\realmult {\fatness} {6.25} \yrad
		\move (0 0) \lellip rx:{\xrad} ry:{\yrad}
		\move (40 -15) \lellip rx:{\xrad} ry:{\yrad}
		\move (40 15) \lellip rx:{\xrad} ry:{\yrad}
		\realadd { 0} { \yrad} \outB
		\realadd { 0} {-\yrad} \outA
		\realadd {15} { \yrad} \inD
		\realadd {15} {-\yrad} \inC
		\realadd {-15} { \yrad} \inB
		\realadd {-15} {-\yrad} \inA
		\move (40 {\inD}) \clvec (15 {\inD})(20 {\outB})(0 {\outB}) 
		\move (40 {\inB}) \clvec (22 {\inB})(22 {\inC})(40 {\inC}) 
		\move (40 {\inA}) \clvec (15 {\inA})(20 {\outA})(0 {\outA}) 
		\savepos (40 -15)(*ex *ey)
	\esegment
	\move (*ex *ey)
}

\newcommand{\handleFA}{%
       \realmult {\fatness} {2.5} \xrad
        \realmult {\fatness} {6.25} \yrad
	\bsegment
		\rmove (40 -15)
    \bsegment
         \move (-40 15) \lellip rx:{\xrad} ry:{\yrad}
        \realadd {15} { \yrad} \outB
        \realadd {15} {-\yrad} \outA
        \realadd {30} { \yrad} \inD
        \realadd {30} {-\yrad} \inC
        \realadd { 0} { \yrad} \inB
        \realadd { 0} {-\yrad} \inA
        \move (0 {\inD}) \clvec (-25 {\inD})(-20 {\outB})(-40 {\outB}) 
        \move (0 {\inB}) \clvec (-18 {\inB})(-18 {\inC})(0 {\inC}) 
        \move (0 {\inA}) \clvec (-25 {\inA})(-20 {\outA})(-40 {\outA}) 
        \savepos (0 0)(*ex *ey)
	  \esegment
	  \move (*ex *ey)

    \bsegment
        \move (40 15) \lellip rx:{\xrad} ry:{\yrad}
        \realadd {15} { \yrad} \outB
        \realadd {15} {-\yrad} \outA
        \realadd {30} { \yrad} \inD
        \realadd {30} {-\yrad} \inC
        \realadd { 0} { \yrad} \inB
        \realadd { 0} {-\yrad} \inA
        \move (0 {\inD}) \clvec (25 {\inD})(20 {\outB})(40 {\outB}) 
        \move (0 {\inB}) \clvec (18 {\inB})(18 {\inC})(0 {\inC}) 
        \move (0 {\inA}) \clvec (25 {\inA})(20 {\outA})(40 {\outA}) 
        \savepos (40 15)(*ex *ey)
    \esegment
    \move (*ex *ey)

	\esegment
    \move (*ex *ey)
}

\newcommand{\intextmultFA}{%
\mathop{
  \raisebox{-2.2pt}{\begin{texdraw}
  \drawdim pt \linewd 0.3 \setunitscale 0.3
  \multFA
  \move (-6 0)
  \move (46 0) 
  \end{texdraw}}
}
}

\newcommand{\intextcomultFA}{%
  \raisebox{-2.3pt}{\begin{texdraw}
  \drawdim pt \linewd 0.3 \setunitscale 0.3
  \comultFA
  \move (-6 0)
  \move (46 0) 
  \end{texdraw}}
}
\newcommand{\intextunitFA}{%
\mathop{
  \raisebox{0.9pt}
  {\begin{texdraw}
  \drawdim pt \linewd 0.3 \setunitscale 0.45
  \unitFA
  \move (6 0)
  \end{texdraw}}
}
}
\newcommand{\intextcounitFA}{%
\mathop{
  \raisebox{0.9pt}
  {\begin{texdraw}
  \drawdim pt \linewd 0.3 \setunitscale 0.45
  \counitFA
  \move (-6 0)
  \end{texdraw}}
}
}
\newcommand{\intextidFA}{%
\mathop{
  \raisebox{1pt}
  {\begin{texdraw}
  \drawdim pt \linewd 0.3 \setunitscale 0.35
  \idFA
  \move (6 0)
  \end{texdraw}}
}
}
\newcommand{\intexthandleFA}{%
\mathop{
  \raisebox{-2.3pt}
  {\begin{texdraw}
  \drawdim pt \linewd 0.3 \setunitscale 0.3
  \handleFA
  \move (-6 0)
  \move (86 0)
  \end{texdraw}}
  }
}

\newcommand{\SFA}{
	\bsegment
		\move ( 0  0) \lellip rx:3 ry:7.5
		\move ( 0 30) \lellip rx:3 ry:7.5
		\move (40 45) \lellip rx:3 ry:7.5
		\move (40 15) \lellip rx:3 ry:7.5
		\move (40 -15) \lellip rx:3 ry:7.5
		\move (80  0) \lellip rx:3 ry:7.5
		\move (80 30) \lellip rx:3 ry:7.5

    \move (0 37.5) \clvec (20 37.5)(15 52.5)(40 52.5)
                   \clvec (65 52.5)(60 37.5)(80 37.5) 
    \move (0 -7.5) \clvec (20 -7.5)(15 -22.5)(40 -22.5)
                   \clvec (65 -22.5)(60 -7.5)(80 -7.5) 
		\move (80 22.5) \clvec (55 22.5)(60 37.5)(40 37.5)
                    \clvec (22 37.5)(22 22.5)(40 22.5)
                    \clvec (65 22.5)(60 7.5)(80 7.5) 
		\move (0 7.5) \clvec (25 7.5)(20 -7.5)(40 -7.5)
                    \clvec (58 -7.5)(58 7.5)(40 7.5)
                    \clvec (15 7.5)(20 22.5)(0 22.5) 
		\savepos (80 0)(*ex *ey)
	\esegment
  \move (*ex *ey)
}

\newcommand{\onedot}[1]{
  \bsegment
	\move (0 0) \fcir f:0 r:2
	\esegment
}

\newcommand{\boundaryFA}{%
  \bsegment
    \realmult {\fatness} {2.5} \xrad
    \realmult {\fatness} {6.25} \yrad
    \move (0 0) \lellip rx:{\xrad} ry:{\yrad}
    \savepos (0 0)(*ex *ey)
  \esegment
  \move (*ex *ey)
}

\newcommand{\multNB}{%
	\bsegment
		\realmult {\fatness} {2.5} \xrad
		\realmult {\fatness} {6.25} \yrad
		\move (0 0) 
		\move (0 30) 
		\move (40 15) 
		\realadd {15} { \yrad} \outB
		\realadd {15} {-\yrad} \outA
		\realadd {30} { \yrad} \inD
		\realadd {30} {-\yrad} \inC
		\realadd { 0} { \yrad} \inB
		\realadd { 0} {-\yrad} \inA
		\move (0 {\inD}) \clvec (25 {\inD})(20 {\outB})(40 {\outB}) 
		\move (0 {\inB}) \clvec (18 {\inB})(18 {\inC})(0 {\inC}) 
		\move (0 {\inA}) \clvec (25 {\inA})(20 {\outA})(40 {\outA}) 
		\savepos (40 15)(*ex *ey)
	\esegment
	\move (*ex *ey)
}

\newcommand{\comultNB}{%
	\bsegment
		\realmult {\fatness} {2.5} \xrad
		\realmult {\fatness} {6.25} \yrad
		\move (0 0) 
		\move (40 -15) 
		\move (40 15) 
		\realadd { 0} { \yrad} \outB
		\realadd { 0} {-\yrad} \outA
		\realadd {15} { \yrad} \inD
		\realadd {15} {-\yrad} \inC
		\realadd {-15} { \yrad} \inB
		\realadd {-15} {-\yrad} \inA
		\move (40 {\inD}) \clvec (15 {\inD})(20 {\outB})(0 {\outB}) 
		\move (40 {\inB}) \clvec (22 {\inB})(22 {\inC})(40 {\inC}) 
		\move (40 {\inA}) \clvec (15 {\inA})(20 {\outA})(0 {\outA}) 
		\savepos (40 -15)(*ex *ey)
	\esegment
	\move (*ex *ey)
}

\newdimen\tableauside\tableauside=1.0ex
\newdimen\tableaurule\tableaurule=0.4pt
\newdimen\tableaustep
\def\phantomhrule#1{\hbox{\vbox to0pt{\hrule height\tableaurule width#1\vss}}}
\def\phantomvrule#1{\vbox{\hbox to0pt{\vrule width\tableaurule height#1\hss}}}
\def\sqr{\vbox{%
  \phantomhrule\tableaustep
  \hbox{\phantomvrule\tableaustep\kern\tableaustep\phantomvrule\tableaustep}%
  \hbox{\vbox{\phantomhrule\tableauside}\kern-\tableaurule}}}
\def\squares#1{\hbox{\count0=#1\noindent\loop\sqr
  \advance\count0 by-1 \ifnum\count0>0\repeat}}
\def\tableau#1{\vcenter{\offinterlineskip
  \tableaustep=\tableauside\advance\tableaustep by-\tableaurule
  \kern\normallineskip\hbox
    {\kern\normallineskip\vbox
      {\gettableau#1 0 }%
     \kern\normallineskip\kern\tableaurule}%
  \kern\normallineskip\kern\tableaurule}}
\def\gettableau#1 {\ifnum#1=0\let\next=\null\else
  \squares{#1}\let\next=\gettableau\fi\next}

\usepackage{verbatim}
\usepackage{amsmath,amsthm,amsfonts}

\newcommand{\cnums} {{\mathbb C}}          
\newcommand{\znums} {{\mathbb Z}}		
\newcommand{\qnums} {{\mathbb Q}}		

\newcommand{\Tr}{\operatorname{tr}}
\newcommand{\tr}{\operatorname{tr}}

\renewcommand{\O}{\mathcal{O}}
\renewcommand{\P}{\mathbb{P}}
\newcommand{\M}{\overline{M}^{\bullet}_{h}}
\newcommand{\E}{\mathbb{E}}

\newcommand{\sss}[1]{\scriptscriptstyle{#1}}

\newcommand{\combinatfactor}{\mathfrak{z}}

\newcommand{\vertline}{\operatorname{|}}

\newcommand{\vac}{v_\emptyset}

\DeclareMathOperator{\Hilb}{Hilb}
\DeclareMathOperator{\Sym}{Sym}
\newcommand{\zz}{{\combinatfactor}}
\newcommand{\cF}{\mathcal{F}}
\newcommand{\Q}{\mathbb{Q}}
\newcommand{\Z}{\mathbb{Z}}
\newcommand{\lv}{\left |}

\newcommand{\lang}{\left\langle}
\newcommand{\rang}{\right\rangle}
\newcommand{\MM}{\mathsf{M}}

\newcommand{\TT}{T^\pm}

\newtheorem{thm}{Theorem}[section]
\newtheorem{theorem}[thm]{Theorem}
\newtheorem{cor}[thm]{Corollary}

\newtheorem{lemma}[thm]{Lemma}

\newtheorem{proposition}[thm]{Proposition}

\newtheorem{definition}[thm]{Definition}

\newtheorem{assumption}{Assumption}
\newtheorem{remark}{Remark}

\begin{document}

\maketitle 

\begin{abstract}
The local Gromov-Witten theory of curves is solved by localization and
degeneration methods. Localization is used for the exact evaluation of
basic integrals in the local Gromov-Witten theory of $\P^1$. A TQFT
formalism is defined via degeneration to capture higher genus
curves. Together, the results provide a compete and effective
solution.

The local Gromov-Witten theory of curves is equivalent to the local
Donaldson-Thomas theory of curves, the quantum cohomology of the
Hilbert scheme points of $\cnums^2$, and the orbifold quantum
cohomology the symmetric product of $\cnums^2$.  The results of the
paper provide the local Gromov-Witten calculations required for the
proofs of these equivalences.
\end{abstract}

\tableofcontents

\section{Introduction}

\subsection{Local Gromov-Witten theory}

The Gromov-Witten theory of threefolds, particularly Calabi-Yau threefolds,
is a very rich subject.  The study of \emph{local
theories}, Gromov-Witten theories of non-compact targets, has
revealed much of the structure.  Let $X$ be a complete,
nonsingular, irreducible
curve of genus $g$ over $\cnums$, and let $$N\to X$$ be a rank 2
vector bundle with $\det N\cong K_{X}$. Then $N$ is a non-compact
Calabi-Yau\footnote{We call any quasi-projective
threefold with trivial canonical bundle Calabi-Yau.} threefold, and the
Gromov-Witten theory, defined and studied in 
\cite{Br-Pa-rigidity, Br-Pa,Br-Pa-TQFT,
Fa-Pa,Pandharipande-degenerate-contributions}, is
called the \emph{local Calabi-Yau theory of $X$}.  We study here the local
theory of curves without imposing the Calabi-Yau condition $\det N\cong
K_{X}$ on the bundle $N$.

The study of non Calabi-Yau local theories has several advantages.  The
calculations of \cite{Fa-Pa,Pandharipande-degenerate-contributions,
Pandharipande-ICM} predict a uniform structure for all threefold theories
closely related to the Calabi-Yau case.  The introduction of non Calabi-Yau
bundles $N$ yields a more flexible mathematical framework in which new
methods arise.  We present a complete solution of the local Gromov-Witten
theory of curves.
The result requires a nonsingularity statement
proven in the Appendix with C. Faber and A. Okounkov.

The space of curves in a 
Calabi-Yau threefold $Y$ is always of virtual
dimension 0. After  suitable (and certainly
non-algebraic) deformation of the geometry
of $Y$, 
we may expect to
find only isolated curves and their multiple covers ---
though no complete statement has yet been proven.

The Gromov-Witten theory of $Y$ may then be viewed as
an enumeration of the isolated curves {\em together} with
a Gromov-Witten theory of local type
for the multiple covers. When defined, the latter theory
should be closely related to the local Gromov-Witten theory
of curves studied here  
\cite{Br-Pa-rigidity}.

\subsection{Equivalences}
The local Gromov-Witten theory of curves is of substantial interest
beyond the original motivations. The local theory may
be viewed as an exactly solved quantum deformation of the 
Hurwitz question of enumerating
ramified coverings of curves. In fact, the
solution has been discovered to arise in many different
geometry contexts.

Our study  of the local Gromov-Witten theory of curves
is a starting point for
several lines of inquiry:
\begin{enumerate}
\item[(i)]
The Gromov-Witten/Donaldson-Thomas correspondence of \cite{MNOP1,MNOP2} may
be naturally studied in the context of local theories.  
Our results together with \cite{dtlc} prove the correspondence
for local theories of curves, see Section 
\ref{subsec: GW/DT for local curves}.

\item[(ii)]
The local theory of the {\em trivial} rank 2 bundle over $\P^1$ 
is equivalent to the quantum cohomologies of the Hilbert
scheme $\Hilb^n(\cnums^2)$ and 
the orbifold
$(\cnums^2)^{n}/S_n$.
Our results here together with
 \cite{Bryan-Graber, Ok-Pan-Hilb} prove the equivalences,
see Section \ref{sec: further directions}.
\end{enumerate}
We expect further connections will likely be found in the future.

\subsection{Results}
Let $N$ be a rank 2 bundle on a curve $X$ of genus $g$.
We assume $N$ is decomposable as a direct sum of
line bundles,
\begin{equation}\label{ffg}
N=L_1\oplus L_2.
\end{equation}
The splitting determines a scaling action of a 2-dimensional torus 
$$T=\cnums^* \times \cnums^*$$ on $N$.
The {\em level} of the splitting is the pair of integers
$(k_1,k_2)$ where,
$$k_i= {\text {deg}}(L_i).$$
Of course, the scaling action and the level 
depend upon the splitting \eqref{ffg}.

The Gromov-Witten residue invariants of $N$, defined in
Section~\ref{subsec: GW residue invariants of N}, take values in the
localized equivariant cohomology ring of $T$ generated by $t_1$ and
$t_2$.  The basic objects of study in our paper are the
\emph{partition functions}
\[
{\mathsf{GW}}_d(g\vertline  k_{1},k_{2}) \in \qnums (t_{1},t_{2})((u)),
\]
the generating functions for the degree $d$ residue invariants of
$N$. Here, $u$ parameterizes the domain genus.

The residue invariants specialize to the local invariants of $X$ in a
Calabi-Yau threefold defined in \cite{Br-Pa,Br-Pa-TQFT} if the level satisfies
\[
k_{1}+k_{2}=2g-2
\]
and the variables are equated,
\[
t_1=t_2.
\] 
Equating the variables is equivalent to considering the residue theory
of $N$ with respect to the diagonal action of a 1-dimensional
torus.

For the Gromov-Witten residue invariants of $N$, we develop a gluing theory
in Section \ref{ctt} 
following \cite{Br-Pa-TQFT}. The interpretation of the local
theory as TQFT is discussed in Section \ref{tqft}.
In Sections \ref{ss} - \ref{ctt} and the Appendix, the gluing relations,
together with a few basic integrals, are proven to
determine the full local theory of curves.
The level freedom of the theory
plays an essential role.  We provide explicit formulas in Sections
\ref{sec: s1=-s2 limit} and \ref{sec: explicit formulas}.

A parallel equivariant Donaldson-Thomas residue theory can be defined for
the threefold $N$.  We conjecture a
Gromov-Witten/Donaldson-Thomas correspondence
for equivariant residues in the framework of
\cite{MNOP1,MNOP2}, see Section~\ref{sec: GW/DT correspondence for
residues}. An important consequence of our theory is Theorem 
\ref{thm: Z is a rational function of q}. After suitable normalization,
${\mathsf {GW}}_d(g\vertline  k_{1},k_{2})$ is a
\emph{rational} function of the variables $t_1$, $t_2$, and
\[
q=-e^{iu}.
\]
The result verifies a
 prediction of the GW/DT correspondence, see Conjecture~2R of
Section~\ref{sec: GW/DT correspondence for residues}.

The residue invariants of $N$ are of special interest when the variable
reduction,
\[
t_1+t_2=0,
\]
is taken. The reduction is equivalent to considering the residue theory of
$N$ with respect to the {\em anti-diagonal} action of a 1-dimensional
torus.  In Theorem \ref{thm: s1=-s2 general formula}, 
we obtain a general closed formula for the partition
function in the anti-diagonal case.

If we additionally specialize to the Calabi-Yau case, our formula is
particularly attractive. The residue partition function here is simply a
$Q$-deformation of the classical formula for unramified covers (see
Corollary~\ref{cor: s1=-s2 limit, CY case}):
\[
{\mathsf{GW}}_d 
(g\vertline  k,2g-2-k)= (-1)^{d (g-1-k)}\sum _{\rho }\left(\frac{d!}{\dim
_{Q}\rho } \right)^{2g-2} Q^{-c_{\rho } (g-1-k)}
\]
where $Q=e^{iu}$ and the sum is over partitions. With the
anti-diagonal action, $N$ is \emph{equivariantly} Calabi-Yau.

Using the above formula, Aganagic, Ooguri, Saulina, and Vafa have
recently found that the local Gromov-Witten theory of curves is
closely related to $q$-deformed 2D Yang-Mills theory and bound states
of BPS black holes \cite{Aganagic-Ooguri-Saulina-Vafa,Vafa-04-2dYang-Mills}.

The anti-diagonal action is exactly {\em opposite} to the original
motivations of the project. It would be very interesting to find
connections between the anti-diagonal case and the original questions of
the Gromov-Witten theory of curves in Calabi-Yau threefolds.

\subsection{Acknowledgments} The authors thank G. Farkas, T. Graber,
A. Greenspoon, S. Katz, J. Kock, C. Teleman, M. Thaddeus, C. Vafa, and
R. Vakil for valuable discussions. We thank J. Kock for the use of his
cobordism \LaTeX \@ macros. A first draft of the Appendix was
completed during a visit by C. Faber to Princeton in the summer of
2004. 

J.~B. was partially supported by the NSERC, the
Clay Institute, and the Aspen Institute.  R.~P. was partially supported by
the Packard foundation and the NSF.

\section{The residue theory}

\subsection{Gromov-Witten residue invariants} 

Let $Y$ be a nonsingular, {\em quasi-projective}, algebraic threefold. 
Let $\M  (Y,\beta )$ denote the moduli space
of stable maps 
\[
f:C\to Y
\]
of genus $h$ and degree $\beta\in H_{2} (Y,\znums )$.
The superscript $\bullet$ indicates the
possibility of disconnected domains $C$. 
We require $f$ to be nonconstant on each connected component of $C$.
The genus, $h(C)$, is defined by  
\[
h (C)= 1-\chi (\O _{C})
\]
and may be negative.

Let $Y$ be equipped with an action by an algebraic torus $T$.
We will define Gromov-Witten residue invariants under the following
assumption.
\begin{assumption}\label{assuption: T fixed locus of GW space is compact}
The $T$-fixed point set $\M (Y,\beta )^{T}$ is compact.
\end{assumption}

We motivate the definition of the residue invariants of $Y$
as follows. We would like to define the reduced Gromov-Witten partition
function ${\mathsf Z}' (Y)_\beta$ as a generating function of the
integrals of the identity class over the moduli spaces of maps,
\begin{equation}\label{eqn: moral defn of Zgw}
{\mathsf Z}' (Y)_\beta\text{ ``}=\text{'' }\sum _{h\in \znums
}u^{2h-2}\int _{[\M (Y,\beta ) ]^{vir}}1 .
\end{equation}
However, the integral on the right might not be well-defined if $Y$ is
not compact. 

If $Y$ has trivial canonical bundle and $\M (Y,\beta )$ is compact,
then the integral \eqref{eqn: moral defn of Zgw} is well-defined. The
resulting series ${\mathsf Z}' (Y)_\beta$ is then the usual
reduced partition function for the degree $\beta$ disconnected
Gromov-Witten invariants\footnote{We follow the notation of
\cite{MNOP1,MNOP2} for the reduced partition function.  The prime
indicates the removal of the degree 0 contributions.  In
\cite{MNOP1,MNOP2}, the moduli space $\M (Y,\beta )$ is denoted by
$\overline{M}'_{h} (Y,\beta )$. However, to maintain notational
consistency with \cite{Br-Pa-TQFT}, we will {\em not} adopt the latter
convention.} of $Y$. We can use the virtual localization formula to
express ${\mathsf Z}'(Y)_\beta$ as a residue integral over the
$T$-fixed point locus.

More generally, under Assumption~\ref{assuption: T fixed locus of GW space
is compact}, the series ${\mathsf Z}'(Y)_\beta$ can be {\em defined}
via localization.

\begin{definition}\label{defn: Zgw'}
The reduced partition function for the degree $\beta $ residue
Gromov-Witten invariants of $Y$ is defined by:
\begin{equation}\label{eqn: definition of Zgw}
{\mathsf Z}'(Y)_\beta =\sum _{h\in \znums }u^{2h-2}\int
_{[\M (Y,\beta )^{T}]^{vir}}\frac{1}{e (\operatorname{Norm}^{vir})}.
\end{equation}
\end{definition}

The $T$-fixed
part of the perfect obstruction theory for $\M (Y,\beta )$ induces a
perfect obstruction theory for $\M (Y,\beta )^{T}$ and hence a virtual
class \cite{Gr-Pa}. 
The 
 equivariant virtual normal bundle
of the embedding, $$\M (Y,\beta )^{T}\subset \M (Y,\beta ), $$ is
 $\operatorname{Norm}^{vir}$  with
equivariant Euler class
$e(\operatorname{Norm}^{vir})$. The integral in \eqref{eqn: definition of Zgw}
denotes equivariant push-forward to a point.

Let $r$ be the rank of $T$, and 
let $t_{1}, \ldots, t_r$ be generators for the equivariant cohomology
of $T$,
\[
H^{*}_{T} (\text{pt})\cong \qnums [t_{1},\ldots,t _{r}].
\]
By Definition~\ref{defn: Zgw'}, $ {\mathsf Z}'(Y)_\beta $ is a Laurent
series in $u$ with coefficients given by rational functions of the
variables $t_{1},\ldots,t_{r}$ of homogeneous degree equal to minus the
virtual dimension of $\M(Y,\beta)$.

\subsection{Gromov-Witten residue invariants of $N$} \label{subsec: GW residue invariants of N} 

Let $X$ be a nonsingular, irreducible, projective curve of genus $g$.  Let
$$N=L_1\oplus L_2$$ be a rank 2 bundle on $X$.  The residue invariants
of the threefold $N$ with respect to the 2-dimensional scaling torus
action can be written in terms of integrals over the moduli space of
maps to $X$.

The residue theory may be considered for the 1-dimensional scaling
torus action on an \emph{indecomposible} rank 2 bundle $N$. Since every
rank 2 bundle is equivariantly deformation equivalent to a
decomposable bundle, the residue invariants of indecomposable bundles
are specializations of the split case.

A stable map to $N$ which is $T$-invariant must factor through the zero
section. Hence,
\[
\M (N,d[X])^{T}\cong \M (X,d).
\]
Moreover, the $T$-fixed part of the perfect obstruction theory of $\M
(N,d[X])$, restricted to $\M (N,d[X])^{T}$, is exactly the usual perfect
obstruction theory for $\M (X,d)$. Hence,
\[
[\M (N,d[X])^{T}]^{vir}\cong [\M (X,d)]^{vir}.
\]

The virtual normal bundle of $\M (N,d[X])^{T}\subset \M (N,d[X])$,
considered as an element of $K$-theory on $\M (X,d)$, is given by
\[
\operatorname{Norm}^{vir}=R^{\bullet }\pi _{*}f^{*} (L_{1}\oplus L_{2})
\]
where
\[
\begin{diagram}[height=0.8cm]
U&\rTo^{f}&X\\
\dTo^{\pi }&&\\
\M (X,d)&&
\end{diagram}
\]
is the universal diagram for $\M (X,d)$. 

The reduced Gromov-Witten partition function of the residue invariants
may be written in the following form via equivariant integration:
\[
{\mathsf Z}'_d(N)=\sum _{h\in \znums }u^{2h-2}\int
_{[\M (X,d)]^{vir}}e (-R^{\bullet }\pi _{*}f^{*} (L_{1}\oplus L_{2})).
\]
We will be primarily interested in a partition function with a shifted
exponent,
\[
{\mathsf{GW}}_d (g\vertline  k_{1},k_{2})
=u^{d (2-2g+k_{1}+k_{2})}\ {\mathsf Z}'_d(N).
\]
The shift can be interpreted geometrically as
$$\int_{d[X]} c_1(T_N) = d(2-2g+k_1+k_2),$$
where $T_N$ is the tangent bundle of the threefold $N$.

The explicit dependence on the equivariant parameters $t_{1}$
and $t_{2} $ may be written as follows. Let $b_{1}$ and $b_{2}$ 
be non-negative integers 
satisfying
$$
b_{1}+b_{2}=2h-2+d (2-2g)$$
where $2h-2+d(2-2g)$ is the  virtual dimension of $\M(X,d)$. 
Let
\[
{\mathsf{GW}}^{b_{1},b_{2}}_d(g\vertline  k_{1},k_{2})
=\int _{[\M (X,d)]^{vir}}c_{b_{1}} (-R^{\bullet }\pi _{*}f^{*}L_{1})c_{b_{2}} 
(-R^{\bullet }\pi _{*}f^{*}L_{2}),
\]
where $\int$ here denotes {\em ordinary} integration.
The equivariant Euler class $e (-R^{\bullet }\pi _{*}f^{*} (L_{1}\oplus
L_{2}))$ is easily expressed in terms of the equivariant parameters and the
{\em ordinary} Chern classes of $-R^{\bullet }\pi _{*}f^{*} (L_{1})$ and
$-R^{\bullet }\pi _{*}f^{*} (L_{2})$,
\begin{multline*}
{\mathsf{GW}}_d (g\vertline  k_{1},k_{2})= \\
u^{d (k_{1}+k_{2})}t_{1}^{d (g-1-k_{1})}t_{2}^{d (g-1-k_{2})} \sum
_{b_{1},b_{2}=0}^{\infty } u^{b_{1}+b_{2}}t_{1}^{\frac{1}{2}(b_{2}-b_{1})}
t_{2}^{\frac{1}{2}(b_{1}-b_{2})}{\mathsf{GW}}^{b_{1},b_{2}}_d
(g\vertline  k_{1},k_{2}).
\end{multline*}
Since $b_{1}+b_{2}$ is even, the exponents of $t_{1}$ and $t_{2}$ are
integers. We see that ${\mathsf{GW}}_d (g\vertline k_{1},k_{2})$ 
is a Laurent
series in $u$ with coefficients given by rational functions of $t_{1}$
and $t_{2}$ of homogeneous degree $d (2g-2-k_{1}-k_{2})$.

\section{Gluing formulas}
\label{gfs}

\subsection{Notation and conventions for partitions}\label{subsub:
conventions for partitions}

By definition, a partition $\lambda $ is a finite
sequence of positive integers
\[
\lambda = (\lambda _{1}\geq \lambda _{2}\geq  \lambda_{3} \geq \dots )
\]
where 
$$|\lambda |=\sum _{i}\lambda _{i}=d.$$ We use the notation $\lambda
\vdash d$ to indicate that $\lambda $ is a partition of $d$.

The number of parts of $\lambda $ is called the \emph{length} of
$\lambda $ and is denoted $l (\lambda )$. Let $m_{i} (\lambda )$ be
the number of times that $i$ occurs in the partition $\lambda $.  We
may write a partition in the format:
\[
\lambda = (1^{m_{1}}2^{m_{2}}3^{m_{3}}\cdots )\ .
\]
The combinatorial factor,
\[
\combinatfactor (\lambda ) =\prod _{i=1}^{\infty }m_{i} (\lambda )!i^{m_{i}
(\lambda )},
\]
arises frequently.

A partition $\lambda $ is uniquely determined by the associated
Ferrers diagram, which is the collection of $d$ boxes located at
$(i,j)$ where $1\leq j\leq \lambda _{i}$. For example
\[
(3,2,2,1,1)= (1^{2}2^{2}3)=\tableau{3 2 2 1 1}.
\]
The \emph{conjugate partition} $\lambda '$ is obtained by reflecting the
Ferrers diagram of $\lambda $ about the $i=j$ line.

In Section~\ref{sec: s1=-s2 limit}, we will require the following standard
quantities.
Given a box in the Ferrers diagram, $\Box \in \lambda $, define the
\emph{content} $c (\Box )$ to be $i-j$, and the \emph{hooklength} $h (\Box
)$ to be $ \lambda_{i} + \lambda_{j} '-i-j+1$. The total content
\[
c_{\lambda }=\sum _{\Box \in \lambda }c (\Box )
\]
and the total hooklength
\[
\sum _{\Box \in \lambda }h (\Box )
\]
satisfy the following identities (page 11 of \cite{Macdonald}):
\begin{equation}\label{eqn: identities for total content and hook length}
\sum _{\Box \in \lambda }h (\Box ) = n (\lambda  )+n (\lambda ')+d,\quad \quad 
\quad c_{\lambda }=n (\lambda' )-n (\lambda ),
\end{equation}
 where
\[
n (\lambda )=\sum _{i=1}^{l (\lambda )} (i-1)\lambda _{i}.
\]

\subsection{Relative invariants}

To formulate our gluing laws for the residue theory of rank 2 bundles on $X$, 
we require relative versions of the
residue invariants.

Motivated by the symplectic theory of A.-M. Li and Y. Ruan \cite{Li-Ruan},
J. Li has developed an algebraic theory of relative stable maps to a pair
$(X,B)$. This theory compactifies the moduli space of maps to $X$ with
prescribed ramification over a non-singular divisor $B\subset X$,
\cite{Li-relative1,Li-relative2}.  Li constructs a moduli space of
relative stable maps together with a virtual fundamental cycle and
proves a gluing formula.

Consider a degeneration of $X$ to $X_{1}\cup_{B} X_{2}$, the union of
$X_{1}$ and $X_{2}$ along a smooth divisor $B$. The gluing formula
expresses the virtual fundamental cycle of the usual stable map moduli
space of $X$ in terms of virtual cycles for relative stable maps of
$(X_{1},B)$ and $(X_{2},B)$. The theory of relative stable maps has
also been pursued in \cite{EGH, Ionel-Parker-Annals2003,Ionel-Parker00}.

In our case, the target is a non-singular curve $X$ of genus $g$,
and the divisor $B$ is a collection of points $x_{1},\dots ,x_{r}\in X$.

\begin{definition}\label{defn: rel stable maps}
Let $(X,x_{1},\dots x_{r})$ be a fixed non-singular genus $g$ curve with
$r$ distinct marked points.  Let $\lambda_1, \ldots, \lambda_r$ be
partitions of $d$.
Let 
\[
\M (X, \lambda _{1},\dots, \lambda_r)
\]
be the moduli space of genus $h$ relative stable maps (in the sense of
 Li)\footnote{For a formal definition of relative stable maps, we refer to
\cite{Li-relative1} Section~4.} with target $(X,x_{1},\dots ,x_{r})$
satisfying the following:
\begin{enumerate}
\item[(i)] The maps have degree $d$.
\item[(ii)] The maps are ramified over $x_{i}$ with ramification type $\lambda
_{i}$.
\item[(iii)] The domain curves are possibly disconnected, but the map is not
degree 0 on any connected component.
\item[(iv)] The domain curves are \emph{not} marked.
\end{enumerate}
The partition $\lambda _{i}\vdash d$ determines a ramification type over
$x_{i}$ by requiring the monodromy of the cover (considered as a conjugacy
class of $S_{d}$) has cycle type $\lambda _{i}$.  
\end{definition}

Our moduli spaces of relative stable maps differ from Li's in a few
minor ways. For a complete discussion, see \cite{Br-Pa-TQFT}.

We define the relative reduced partition function 
via equivariant integration over spaces of relative stable maps:
\[
{\mathsf Z}'(N)_{\lambda^1 \dots \lambda^r} =
\sum _{h\in \znums }
u^{2h-2}
 \int _{[\M (X,\lambda ^{1},\dots ,\lambda ^{r})]^{vir}}e
(-R^{\bullet }\pi _{*}f^{*} (L_{1}\oplus L_{2})).
\]
Again, we will be primarily interested in a shifted generating function,
$$
{\mathsf{GW}}(g\vertline  k_{1},k_{2})_{\lambda ^{1} \dots \lambda ^{r}}= 
u^{d (2-2g+k_{1}+k_{2}-r)+\sum _{i=1}^{r}l (\lambda ^{i})} \ 
{\mathsf Z}'(N)_{\lambda^1 \dots \lambda^r}.$$
Since the degree $d$ is equal to $|\lambda^i|$, the degree
subscript is redundant in the relative theory.

The  
exponent of $u$ in the partition function 
${\mathsf{GW}}_d(g\vertline  k_1,k_2)$ of the non-relative theory
is 
$$2h-2+ \int_{d[X]}c_1(T_N).$$ 
In the relative theory, the $2h-2$ term in the exponent is
replaced with $2h-2+\sum l (\lambda ^{i})$, the negative Euler
characteristic of the \emph{punctured} domain.  The class $c_1(T_N)$ 
is replaced with the dual of the log canonical class of $N$ with
respect to the relative divisors.
The outcome is the modified exponent of $u$ in the partition function 
${\mathsf{GW}}(g\vertline  k_{1},k_{2})_{\lambda ^{1} \dots \lambda ^{r}}.$

As before, we can make the dependence on $t_{1}$ and $t_{2}$ explicit. Let
$$b_{1}+b_{2}=2h-2+d (2-2g)-   \delta,$$
where 
\[
\delta=\sum _{i=1}^{r} (d-l (\lambda ^{i})).
\]
Here,  $b_1+b_2$ equals the virtual dimension
 of $\M (X,\lambda ^{1}\dots \lambda ^{r}).$
Let
\[
{\mathsf{GW}}^{b_{1},b_{2}} 
(g\vertline  k_{1},k_{2})_{\lambda ^{1} \dots \lambda ^{r}}=\int
_{[\M (X,\lambda ^{1}, \dots ,\lambda ^{r})]^{vir}}c_{b_{1}} (-R^{\bullet
}\pi _{*}f^{*}L_{1})c_{b_{2}} (-R^{\bullet }\pi _{*}f^{*}L_{2}).
\]
Then, we have
\begin{multline}\label{eqn: expression for Zdgk1k2_lambdas}
{\mathsf{GW}} (g\vertline  k_{1},k_{2})_{\lambda ^{1}\dots \lambda ^{r}}=
u^{d (k_{1}+k_{2})}t_{1}^{d (g-1-k_{1})}t_{2}^{d (g-1-k_{2})} \\ \cdot \sum
_{b_{1},b_{2}=0}^{\infty } 
u^{b_{1}+b_{2}}t_{1}^{\frac{b_{2}-b_{1}+\delta}{2}}
t_{2}^{\frac{b_{1}-b_{2}+\delta}{2}}
{\mathsf{GW}}^{b_{1},b_{2}} 
(g\vertline  k_{1},k_{2})_{\lambda ^{1}\dots \lambda ^{r}}.
\end{multline}
Since the parity of $b_{1}+b_{2}$ is the same as $\delta$, the
exponents of $t_{1}$ and $t_{2}$ are integers.

The partition function 
${\mathsf{GW}}(g\vertline k_{1},k_{2})_{\lambda ^{1} \dots
\lambda ^{r}}$ is a Laurent series in $u$ with coefficients given by
rational functions in $t_{1}$ and $t_{2}$ of homogeneous degree
\[
d (2g-2-k_{1}-k_{2})+\delta.
\]

In \cite{Br-Pa-TQFT}, the combinatorial factor $\combinatfactor
(\lambda )$ is used to raise the indices for the relative invariants. For the
residue invariants, an additional factor
$(t_{1}t_{2})^{l (\lambda )}$ must be included. We define:
\begin{equation}\label{eqn: index raising formula}
{\mathsf{GW}} 
(g\vertline  k_{1},k_{2})_{\mu ^{1}\dots \mu ^{s}}^{\nu^{1}\dots \nu
^{t}}= 
{\mathsf{GW}} (g\vertline  k_{1},k_{2})_{\mu
^{1}\dots \mu ^{s},\nu ^{1}\dots \nu ^{t}}
\left(\prod _{i=1}^{t} \combinatfactor (\nu ^{i})
(t_{1}t_{2})^{l (\nu ^{i})} \right).
\end{equation}

\subsection{Gluing formulas}
The gluing formulas are determined by the following result.

\begin{theorem}\label{thm: gluing formula}
For splittings
$g=g'+g''$ and $k_{i}=k_{i}'+k_{i}''$,
\[
{\mathsf{GW}} (g\vertline k_{1},k_{2})_{\mu ^{1}\dots \mu ^{s}}^{\nu ^{1}\dots
\nu ^{t}} =\sum _{\lambda \vdash d} {\mathsf{GW}} (g'\vertline
k_{1}',k_{2}')_{\mu ^{1}\dots \mu ^{s}}^{\lambda }
\mathsf{GW} (g''\vertline
k_{1}'',k_{2}'')_{\lambda }^{\nu ^{1}\dots \nu ^{t}}
\]
and
\[
{\mathsf{GW}}
 (g\vertline  k_{1},k_{2})_{\mu ^{1}\dots \mu ^{s}}=\sum _{\lambda \vdash
d}{\mathsf{GW}} (g-1\vertline  k_{1},k_{2})^{\lambda }_{\mu ^{1}\dots \mu ^{s}\lambda }.
\]
\end{theorem}

\begin{proof}
The proof follows the derivation of the gluing formulas in \cite{Br-Pa-TQFT}. 
The only difference
is the modified metric term
$$\combinatfactor(\lambda) (t_1t_2)^{l(\lambda)}.$$ The first factor,
$\combinatfactor(\lambda)$, is obtained from the degeneration formula
for the virtual class \cite{Li-relative2} as in
\cite{Okounkov-Pandharipande-completed-cycles}.

The second factor, $(t_1t_2)^{l(\lambda)}$,
arises from normalization sequences associated to the fractured domains.
Let 
$$f:C \rightarrow X$$
be an element of $\M (X,\mu ^{1},\dots ,\mu ^{s}, \nu^1,\dots,\nu^t)$.
Consider a reducible degeneration of the target, 
$$X=X'\cup X'',$$
over
which the line bundles $L_1$ and $L_2$ extend with degree splittings
$$k_1=k_1'+k_1'',$$
$$k_2=k_2'+k_2''.$$ 
In a degeneration of type $\lambda$,
the domain curve degenerates,
$$C= C' \cup C'',$$
into components lying over $X'$ and $X''$ and satisfying
$$|C' \cap C''| = l(\lambda).$$
For each line bundle $L_i$, we have a normalization sequence, 
\begin{equation}
\label{rezz} 0 \rightarrow f^*(L_i)|_{C} \rightarrow f^*(L_i)|_{C'}
\oplus f^*(L_i)|_{C''} \rightarrow f^*(L_i)|_{C'\cap C''} \rightarrow
0.
\end{equation}
The last term yields a trivial bundle of rank
$l(\lambda)$
with scalar torus action over the moduli space
of maps of degenerations of type $\lambda$.
The factor $(t_1t_2)^{l(\lambda)}$ is obtained from the
higher direct images of the normalization sequences \eqref{rezz}.
The analysis for irreducible degenerations of $X$ is identical.

The exponent of $u$ in the series 
${\mathsf{GW}} 
(g\vertline  k_{1},k_{2})_{\mu ^{1}\dots \mu ^{s}}^{\nu ^{1}\dots \nu
^{t}}$ of relative
invariants has been precisely chosen to respect the
gluing rules.
\end{proof}

\section{TQFT formulation of gluing laws} 
\label{tqft}
\subsection{Overview}
The gluing structure of the residue theory of rank 2 bundles on curves
is most concisely formulated as a functor of tensor categories,
\[
\mathbf{{GW}}(-):2\mathbf{Cob}^{L_{1},L_{2}}\to R\mathbf{mod}.
\]

The deformation invariance of the residue theory allows for a
topological formulation of the gluing structure.

Our discussion follows Sections 2 and 4 of \cite{Br-Pa-TQFT} and draws from
Chapter~1 of \cite{Kock:FA-2DTQFT}.  Modifications of the categories have
to be made to accommodate the more complicated objects studied here.

\subsection{ $2\mathbf{Cob}$ and $2\mathbf{Cob}^{L_{1},L_{2}}$}

We first define the category $2\mathbf{Cob}$
of 2-cobordisms. The objects of 
$2\mathbf{Cob}$
 are compact oriented
1-manifolds, or equivalently, finite unions of oriented circles.  Let
$Y_1$ and $Y_2$ be objects of the category.  A morphism,
\[
Y_{1}\rightarrow Y_{2},
\] 
is an equivalence class of oriented cobordisms $W$ from $Y_{1}$ to
$Y_{2}$. Two cobordisms are equivalent if they are diffeomorphic by a
boundary preserving oriented diffeomorphism. Composition of morphisms
is obtained by concatenation of the corresponding cobordisms. The
tensor structure on the category is given by disjoint union.

The category $2\mathbf{Cob}^{L_{1},L_{2}}$ is defined to have the same objects as
$2\mathbf{Cob}$. A morphism in $2\mathbf{Cob}^{L_{1},L_{2}},$ 
$$Y_{1} \rightarrow Y_{2},$$ 
is an equivalence class of triples $(W,L_{1},L_{2})$ where $W$
is an oriented cobordism from $Y_{1}$ to $Y_{2}$ and $L_{1},L_{2}$ are
complex line bundles on $W$, trivialized on $\partial W$. The triples
$(W,L_{1},L_{2})$
and  $(W',L'_{1},L'_{2})$ are equivalent if there exists a boundary preserving
oriented diffeomorphism, $$f:W\to W',$$ and bundle isomorphisms $$L_{i}\cong
f^{*}L'_{i}.$$ Composition is given by concatenation of the cobordisms and
gluing of the bundles along the concatenation using the trivializations.

The isomorphism class of $L_{i}$ is determined by the Euler class $$e
(L_{i})\in H^{2} (W,\partial W),$$ which assigns an integer to each
component of $W$. For a connected cobordism $W$, we refer to the pair of
integers $(k_{1},k_{2})$, determined by the Euler classes of $L_{1}$ and
$L_{2}$, as the \emph{level}. Under concatenation, the
levels simply add. For example:

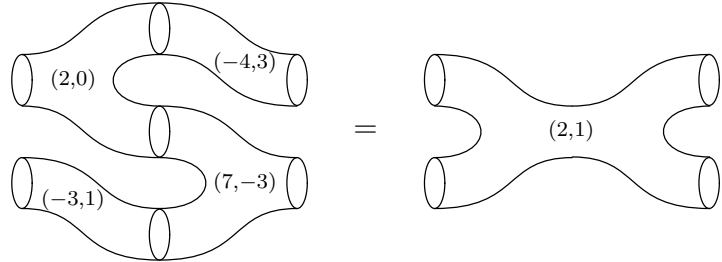
\begin{center}
\begin{texdraw}\setunitscale 1.3
\SFA 
\move (120 0)
\boundaryFA \multNB \comultNB
\move (120 30)\boundaryFA 
\move (200 0)\boundaryFA 
\move (200 30)\boundaryFA 
\htext (100 15) {$=$}
\htext (15 -5) {$\scriptstyle{(-3,1)}$}
\htext (15 30) {$\scriptstyle{(2,0)}$}
\htext (65 0) {$\scriptstyle{(7,-3)}$}
\htext (65 35) {$\scriptstyle{(-4,3)}$}
\htext (160 15) {$\scriptstyle{(2,1)}$}
\end{texdraw}
\end{center}

The empty manifold is a distinguished object in $2\mathbf{Cob}$
and $2\mathbf{Cob}^{L_{1},L_{2}}$. A morphism in
$2\mathbf{Cob}^{L_{1},L_{2}}$ from the empty manifold to itself is given by
a compact, oriented, closed 2-manifold $X$ together with a pair of complex
line bundles $L_{1}\oplus L_{2}\to X$.

The full subcategory of $2\mathbf{Cob}^{L_{1},L_{2}}$ obtained by
restricting to level $(0,0)$ line bundles is clearly isomorphic to the
category $2\mathbf{Cob}$.

More generally, we obtain an embedding $2\mathbf{Cob}\subset
2\mathbf{Cob}^{L_{1},L_{2}}$ for any fixed integers $(a,b)$ by requiring
the level of any connected cobordism to be $(a\chi ,b\chi )$ where $\chi $
is the Euler characteristic of the cobordism.

If $a+b=-1$, such an embedding is termed {\em Calabi-Yau}
since the threefold $$L_{1}\oplus L_{2}\to X$$
has numerically trivial canonical class if $$\deg (L_{1})+ \deg (L_{2})=-\chi.$$

\subsection{Generators for $2\mathbf{Cob}$ and
$2\mathbf{Cob}^{L_{1},L_{2}}$} 

The category $2\mathbf{Cob}$ is generated by the morphisms
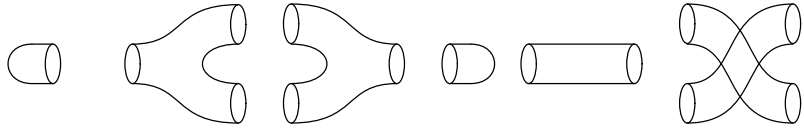
\[
\begin{texdraw}\setunitscale 1.0
\unitFA
\move (30 0)
\comultFA 
\move (90 -15)
\multFA 
\move (150 0)
\counitFA 
\move (180 0)
\idFA 
\move (240 -15)
\twistFA 
\end{texdraw}
\]

In other words, any morphism (cobordism) can be obtained by taking
compositions and tensor products (concatenations and disjoint unions) of
the above list (Proposition~1.4.13 of \cite{Kock:FA-2DTQFT}).

The category $2\mathbf{Cob}^{L_{1},L_{2}}$ is then clearly generated by the
morphisms
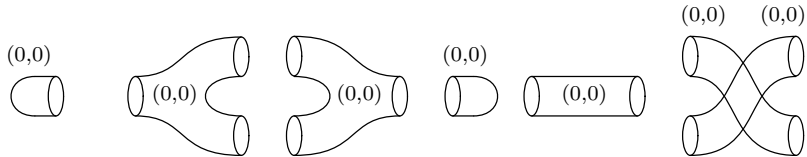
\begin{center}
\begin{texdraw}\setunitscale 1.0
\unitFA
\move (30 0)
\comultFA 
\move (90 -15)
\multFA 
\move (150 0)
\counitFA 
\move (180 0)
\idFA 
\move (240 -15)
\twistFA 
\htext (-10 15) {$\scriptstyle{(0,0)}$}
\htext (45 0) {$\scriptstyle{(0,0)}$}
\htext (115 0) {$\scriptstyle{(0,0)}$}
\htext (155 15) {$\scriptstyle{(0,0)}$}
\htext (200 0) {$\scriptstyle{(0,0)}$}
\htext (245 30) {$\scriptstyle{(0,0)}$}
\htext (275 30) {$\scriptstyle{(0,0)}$}
\end{texdraw}
\end{center}
along with the morphisms
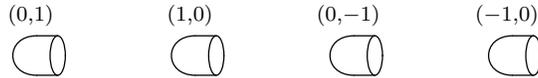
\begin{center}
\begin{texdraw}\setunitscale 1.0
\unitFA 
\move (60 0)
\unitFA 
\move (120 0)
\unitFA 
\move (180 0)
\unitFA 
\htext (-10 15) {$\scriptstyle{(0,1)}$}
\htext (50 15) {$\scriptstyle{(1,0)}$}
\htext (110 15) {$\scriptstyle{(0,-1)}$}
\htext (170 15) {$\scriptstyle{(-1,0)}$}
\end{texdraw}
\end{center}

Let $R$ be a commutative ring with unit, and let $R\mathbf{mod}$ be
the tensor category of $R$-modules.  By a well-known result (see
Theorem 3.3.2 of \cite{Kock:FA-2DTQFT}), a 1+1 dimensional $R$-valued
TQFT, which is by definition a symmetric tensor functor
\begin{equation}
\label{ejf8923}
\mathbf{F}:2\mathbf{Cob}\to R\mathbf{mod},
\end{equation}
is equivalent to a commutative Frobenius algebra over $R$. 

Given a symmetric tensor functor \eqref{ejf8923},
the underlying
$R$-module of the Frobenius algebra is given by
\[
H=\mathbf{F} (S^{1})
\]
and the Frobenius algebra structure is determined as follows:
\[
\begin{array}{lrl}
\text{multiplication}&\mathbf{F} (\intextmultFA ):&H\otimes H\to H\\
\text{unit}&\mathbf{F} (\intextunitFA ):&R\to H\\
\text{comultiplication}&\mathbf{F} (\intextcomultFA ):&H\to H\otimes H\\
\text{counit}&\mathbf{F} (\intextcounitFA ):&H\to R.  
\end{array}
\]

Let $\mathbf{F}$ be a symmetric tensor functor on the larger category
$2\mathbf{Cob}^{L_1,L_2}$,
\[
\mathbf{F}:2\mathbf{Cob}^{L_{1},L_{2}}\to R\mathbf{mod}.
\]
Since the functor $\mathbf{F}$ is determined by the
 values on the generators of $2\mathbf{Cob}^{L_1,L_2}$, the functor 
$\mathbf{F}$
is determined by the level $(0,0)$ Frobenius algebra together with the 
elements

\[
\mathbf{F}\left(\intextunitFA
^{\sss{(0,1)}} \right),\ 
\mathbf{F}\left(\intextunitFA
^{\sss{(1,0)}} \right),\ 
\mathbf{F}\left(\intextunitFA
^{\sss{(0,-1)}} \right),\ 
\mathbf{F}\left(\intextunitFA
^{\sss{(-1,0)}} \right),\ 
\]
Since the latter two elements are the inverses in the Frobenius
algebra of the first two, we obtain half of the following Theorem.

\begin{theorem}\label{thm: Z:2cobL1L2-->Rmod <--> Frob Alg and 2 inv elts}
A symmetric tensor functor 
\[
\mathbf{F}:2\mathbf{Cob}^{L_{1},L_{2}}\to R\mathbf{mod}
\]
is uniquely determined by a commutative Frobenius algebra over $R$ 
for the level $(0,0)$ theory
and two
distinguished, invertible elements 

\[
\mathbf{F}\left(\intextunitFA
^{\sss{(0,-1)}} \right),\ 
\mathbf{F}\left(\intextunitFA
^{\sss{(-1,0)}} \right).
\]
\end{theorem}

\begin{proof}
Uniqueness was proved above. The existence result will not be used in
the paper. We leave the details to 
 the reader.
\end{proof}

\subsection{The functor $\mathbf{GW} (-)$} 

Let $R$ be the ring of Laurent series in $u$ whose coefficients are rational
functions in $s _{1}$ and $s _{2}$,
\[
R=\qnums (t_{1},t_{2})((u)).
\]
The collection of partition functions 
$\mathsf{GW} (g\vertline  k_{1},k_{2})_{\lambda
^{1}\dots \lambda ^{r}}$ of degree $d$ gives rise to a functor
\[
\mathbf{GW}(-):2\mathbf{Cob}^{L_1,L_2}\to R\mathbf{mod}
\]
as follows. Define 
\[
\mathbf{GW} (S^{1})=H=\bigoplus _{\lambda \vdash d}Re_{\lambda }
\]
to be the free $R$-module with basis $\{e_{\lambda } \}_{\lambda \vdash d}$
labeled by partitions of $d$, and let
\[
\mathbf{GW} \left(S^{1} \coprod\dots \coprod S^{1} \right)
=H\otimes \dots \otimes H.
\]

Let $W_{s}^{t} (g\vertline  k_{1},k_{2})$ be the connected genus $g$ cobordism from a
disjoint union of $s$ circles to a disjoint union of $t$ circles, equipped
with line bundles $L_{1}$ and $L_{2}$ of level $(k_{1},k_{2})$. 
We define the $R$-module homomorphism
\[
\mathbf{GW} \left(W_{s}^{t} 
(g\vertline  k_{1},k_{2})\right):H^{\otimes s}\to H^{\otimes t}
\]
by
\[
e_{\eta ^{1}}\otimes \dots \otimes e_{\eta ^{s}}\mapsto \sum _{\mu
^{1}\dots \mu ^{t}\vdash d} {\mathsf{GW}} (g\vertline k_{1},k_{2})_{\eta ^{1}\dots
\eta ^{s}}^{\mu ^{1}\dots \mu ^{t}}e_{\mu ^{1}}\otimes \dots \otimes
e_{\mu ^{t}}.
\]
We extend the definition of $\mathbf{GW} (-)$ to disconnected cobordisms
using tensor products:
\[
\mathbf{GW} \left(W[{1}]\coprod \dots \coprod W[{n}]\right)=\mathbf{GW}
\left(W[{1}])\otimes \dots \otimes \mathbf{GW} (W[{n}]\right).
\]

\begin{theorem}\label{thm: gluing defines TQFT functor}
$\mathbf{GW} (-): 2\mathbf{Cob}^{L_{1},L_{2}}\to
R\mathbf{mod}$ is a well-defined functor.
\end{theorem}
\textsc{Proof:} Following the proof of Proposition~4.1 of
\cite{Br-Pa-TQFT},
the gluing laws imply the following compatibility:
\[
\mathbf{GW} \big((W,L_{1},L_{2})\circ
(W',L'_{1},L'_{2})\big)=\mathbf{GW}(W,L_{1},L_{2})\circ \mathbf{GW}
(W',L'_{1},L'_{2}).
\]
We must also prove that $\mathbf{GW} (-)$ takes identity morphisms to
identity morphisms. Since $W_{1}^{1} (0\vertline 0,0)$ is the identity
morphism from $S^{1}$ to itself in $2\mathbf{Cob}^{L_{1},L_{2}}$, we
require
\begin{equation}\label{jjww}
{\mathsf{GW}} (0\vertline  0,0)_{\mu }^{\nu }=\delta _{\mu}^{\nu }.
\end{equation}
Equation \eqref{jjww} will be proved in Lemma~\ref{lem: charge (0,0) tube}.\qed

\section{Semisimplicity in level $(0,0)$} \label{ss}
\subsection{Rings of definition}
The partition functions for the level $(0,0)$ relative invariants lie
in the ring of power series in $u$,
$$\mathsf{GW}
(g\vertline  0,0)_{\lambda ^{1}\dots\lambda ^{r}} \in \qnums 
(t_{1},t_{2})[[u]],$$
since, by
equation \eqref{eqn: expression for Zdgk1k2_lambdas}, no negative powers of
$u$ appear. 
The level $(0,0)$ relative invariants therefore
determine a commutative
Frobenius algebra over the ring
\[
R=\qnums(t_{1},t_{2})[[u]].
\]

We will require formal square
roots of $t_{1} $ and $t_{2}$. Let $\tilde{R}$ be the complete local ring
of power series in $u$ whose coefficients are rational functions in
$t_{1}^{\frac{1}{2}}$ and $t_{2}^{\frac{1}{2}}$,
\[
\tilde{R}=\qnums(t_{1}^{\frac{1}{2}},t^{\frac{1}{2}}_{2}) [[u]].
\]

\subsection{Semisimplicity} A commutative Frobenius algebra $A$ is
{\em semisimple} if $A$ is isomorphic to a direct sum of 1-dimensional
Frobenius algebras.

\begin{proposition}\label{prop: Zd is semi-simple}
The Frobenius algebra determined by the level $(0,0)$ sector of
$\mathbf{GW} (-)$ in degree $d$ is semisimple over $\tilde{R}$.
\end{proposition}
\textsc{Proof:} 
$\tilde{R}$ is a complete local ring with maximal ideal $m$ generated by $u$. 
Let $F$ be the Frobenius algebra determined by the level $(0,0)$ theory
in degree $d$.
The underlying $\tilde{R}$-module of the Frobenius
algebra $F$,
$$H=\oplus _{\lambda \vdash d}\tilde{R}e_{\lambda },$$
is freely generated.
By Proposition~2.2 of
\cite{Br-Pa-TQFT},  $F$ is semisimple if and only if $F/mF$ is
semisimple over $\tilde{R}/m\tilde{R}$.

The structure constants of the multiplication
in $F/mF$ are given by the $u=0$ specialization of the invariants 
$\mathsf{GW}(0\vertline  0,0)_{\alpha \beta }^{\gamma }$. 
By \eqref{eqn: expression for Zdgk1k2_lambdas}, after the $u=0$ specialization,
only the 
$$b_{1}=b_{2}=0$$
terms remain. The latter are the expected dimension 0 terms with domain genus
\[
2h-2 = d -l (\alpha )-l (\beta )-l (\gamma ).
\]
In the expected dimension 0 case, the moduli space 
$\M (\P ^{1},\alpha, \beta, \gamma )$
is nonsingular of {\em actual} dimension 0. We conclude:
\begin{align*}
\mathsf{GW} (0\vertline  0,0)_{\alpha \beta }^{\gamma }|_{u=0}&= \combinatfactor(\gamma)
(t_{1}t_{2})^{l (\gamma )}\ 
\mathsf{GW} (0\vertline 0,0)_{\alpha \beta \gamma }|_{u=0}\\
&= 
 \combinatfactor(\gamma) 
(t_{1}t_{2})^{\frac{1}{2} (d-l (\alpha )-l (\beta )+l (\gamma ))}
\ \int _{[\M (\P^1,\alpha,\beta ,\gamma)]}1\\
&=\combinatfactor(\gamma)  
(t_{1}t_{2})^{\frac{1}{2} (d-l (\alpha )-l (\beta )+l (\gamma ))}
\ \mathbf{H}_{d}^{\P ^{1}} (\alpha ,\beta ,\gamma ),
\end{align*}
where $\mathbf{H}_{d}^{\P ^{1}} (\alpha ,\beta ,\gamma )$ is the Hurwitz
number of degree $d$ covers of $\P ^{1}$ with prescribed ramification
$\alpha $, $\beta $, and $\gamma $ over the 
points $0,1,\infty \in \P^1$.

Up to factors of $t_1$ and $t_2$, the quotient $F/mF$ is
the Frobenius algebra associated to 
the TQFT studied by Dijkgraaf-Witten and Freed-Quinn
\cite{Dijkgraaf-Witten90,Freed-Quinn}. The latter Frobenius algebra is
 isomorphic to $\qnums [S_{d}]^{S_{d}}$, the center of the group
algebra of the symmetric group, and is well-known to be
semisimple.

We derive below an explicit idempotent basis for $F/mF$ analogous to the
well-known idempotent basis for $\qnums [S_{d}]^{S_{d}}$. 
The formal square roots of $t_{1}$ and $t_{2}$ are required here.

Let $\rho$ be an irreducible representation of $S_d$.  The conjugacy
classes of $S_d$ are indexed by partitions $\lambda$ of size $d$.  Let
$\chi^\rho_\lambda$ denote the trace of $\rho$ on the conjugacy class
$\lambda$.  The Hurwitz numbers are determined by the following
formula:
\[
 \mathbf{H}^{\P ^{1}}_{d} (\alpha ,\beta ,\gamma )=
\sum _{\rho }\frac{d!}{\dim \rho }
\frac{\chi^\rho_\alpha}{\combinatfactor (\alpha )}
\frac{\chi^\rho
_\beta}{\combinatfactor (\beta )}
\frac{\chi^\rho_\gamma }{\combinatfactor (\gamma )},
\]
see, for example, \cite{Okounkov-Pandharipande-completed-cycles} equation {0.8}.
The above sum is over all irreducible representations $\rho$ of $S_d$.

The structure constants for multiplication in $F/mF$ are 
\begin{equation}\label{eqn: u=0 structure constants}
\mathsf{GW} (0\vertline  0,0)^{\gamma }_{\alpha \beta }|_{u=0}
= (t_{1}t_{2})^{\frac{1}{2}
(d-l (\alpha )-l (\beta )+l (\gamma ))} \sum _{\rho }\frac{d!}{\dim \rho
}     \frac{\chi^\rho_\alpha}{\combinatfactor(\alpha)}
\frac{\chi^\rho_\beta}{\combinatfactor
(\beta )}\chi^\rho_\gamma.
\end{equation}
We define a new basis $\{v_{\rho }^{0} \}$ for $F/mF$ by
\begin{equation}\label{eqn: v0 basis}
v_{\rho }^{0}=\frac{\dim \rho }{d!}\sum _{\alpha }
\left(t_{1}^{\frac{1}{2}}t_{2}^{\frac{1}{2}} \right)^{l (\alpha )-d}\chi
^{\rho }_\alpha e_{\alpha }.
\end{equation}
The elements $\{v_{\rho }^{0} \}$ form an idempotent basis:
\[
v_{\rho }^{0}\cdot v_{\rho '}^{0}=\delta _{\rho \rho '}^{\rho ''}v_{\rho
''}^{0}.
\]
By Proposition~2.2 of \cite{Br-Pa-TQFT}, there exists a unique idempotent basis
$\{v_{\rho } \}$ of $F$, such that $v_{\rho }=v_{\rho }^{0} \mod m$. \qed

\begin{remark}
In general, $v_{\rho }\neq v_{\rho }^{0}$ but for the
 anti-diagonal specialization $$t_{1}=-t_{2},$$
the equality $v_{\rho }=v_{\rho }^{0}$ holds (see Section~\ref{sec: s1=-s2
limit}).
\end{remark}

\subsection{Structure}
Semisimplicity leads to a basic structure result.
\begin{theorem}\label{thm: structure theorem}
There exist universal series, 
$\lambda_{\rho },\eta _{\rho
}\in \tilde{R},$ labeled by partitions $\rho$,
for which
\[
\mathsf{GW}_d 
(g\vertline  k_{1},k_{2})
=\sum _{\rho\vdash d} \lambda_{\rho} ^{g-1}\eta_{\rho}
^{-k_{1}}\overline{\eta }_{\rho}^{-k_{2}}.
\]
Here, $\overline{\eta }_{\rho}$ is obtained from $\eta _{\rho}$ by
interchanging $t_{1}$ with $t_{2}$.
\end{theorem}

\textsc{Proof:} Let $\{v_{\rho } \}$ be an idempotent basis for the level
$(0,0)$ Frobenius algebra of $\mathsf{GW} (-)$ in degree $d$.

Define $\lambda _{\rho }$ to be the inverse of the counit evaluated on
$v_{\rho }$:
\[
\lambda _{\rho }^{-1}=\mathsf{GW} \left(\intextcounitFA
^{\sss{(0,0)}} \right) (v_{\rho }).
\]
Equivalently, $\lambda _{\rho }$ is the eigenvalue for the eigenvector
$v_{\rho }$ under the \emph{genus adding operator} $G$:
\[
G=\mathsf{GW} \left(\intexthandleFA^{\sss{(0,0)}} \right):H\to H.
\]

Let $\eta _{\rho }$ (respectively $\overline{\eta }_{\rho }$) be the
coefficient of $v_{\rho }$ in the element
\[
\eta =\mathsf{GW} \left(\intextunitFA^{\sss{(-1,0)}} \right)\in H,\quad
(\text{respectively }\overline{\eta }=\mathsf{GW}
\left(\intextunitFA^{\sss{(0,-1)}} \right)\in H).
\]
Equivalently, $\eta _{\rho }$ (respectively $\overline{\eta }_{\rho }$) is
the eigenvalue for the eigenvector $v_{\rho }$ under the \emph{left 
annihilation operator} (respectively \emph{right annihilation
operator}):
\[
A=\mathsf{GW} \left(\intextidFA^{\sss{(-1,0)}} \right)\quad (\text{respectively
}\overline{A}=\mathsf{GW} \left(\intextidFA^{\sss{(0,-1)}} \right)).
\]
The gluing rules imply: 
\[
{\mathsf{GW}} (g\vertline  k_{1},k_{2})_d=
\tr (G^{g-1}A^{-k_{1}}\overline{A}^{-k_{2}}).
\]
The operators $G$, $A$, and $\overline{A}$ are simultaneously diagonalized
by the basis $\{v_{\rho } \}$, so the Theorem is equivalent to the above
formula.\qed

\section{Solving the theory} \label{ctt}
\subsection{Overview} 
The {\em full local theory of curves} is the set
of all series
\begin{equation}
\mathsf{GW}(g\vertline  k_1,k_2)_{\lambda^1\dots\lambda^r}.
\end{equation}
The functor ${\mathbf {GW}}$ contains the
data of the full local theory.
By Theorem \ref{thm: Z:2cobL1L2-->Rmod <--> Frob Alg and 2 inv elts},
the full local theory is
determined by the following  {\em basic} series:
\begin{equation}
\label{eqn: basic series} \mathsf{GW}(0\vertline  0,0)_{\lambda \mu \nu }\ , \
\mathsf{GW}(0\vertline  0,0)_{\lambda \mu}\ , \ 
\mathsf{GW}(0\vertline  0,0)_\lambda\ ,
\end{equation}
$$ \ \mathsf{GW}(0\vertline  -1,0)_\lambda\ , \ 
 \mathsf{GW}(0\vertline  0,-1)_\lambda.$$

We present a recursive method for calculating the full local
theory of curves using the TQFT formalism.
Four of the basic series,
\[
 \mathsf{GW}(0\vertline  0,0)_{\lambda \mu}\ , 
\ \mathsf{GW}(0\vertline  0,0)_\lambda\ , \
\mathsf{GW}(0\vertline  -1,0)_\lambda\ , \ 
\mathsf{GW}(0\vertline  0,-1)_\lambda\ , 
\]
are determined by closed formulas.  The first two are easily
obtained by dimension considerations (Lemmas \ref{lem: charge (0,0) tube}
and \ref{lem: charge (0,0) cap}).  The last two have been determined in
\cite{Br-Pa-TQFT} in case the equivariant parameters $t_i$ are set to
1. The insertion of the equivariant parameters is straightforward
(Lemma~\ref{lem: the CY cap}).

The level $(0,0)$ pair of pants series
\[
\mathsf{GW}(0\vertline  0,0)_{\lambda \mu  \nu }
\]
are much more subtle. 
The main result of the Appendix (with C. Faber and A. Okounkov) is
the determination of {all} degree $d$ level $(0,0)$ pair of pants series
from the {\em single} series
\begin{equation}\label{nnqq}
\mathsf{GW}(0\vertline  0,0)_{(d),(d),(1^{d-2}2)}
\end{equation}
using the TQFT associativity relations, level $(0,0)$ series of {\em lower
degree}, and Hurwitz numbers of covering genus 0.  A closed formula for
\eqref{nnqq} is derived in Section~\ref{qw223}. The outcome is a
computation of the full local theory of curves via recursions in degree
(Theorem~\ref{thm: reconstruction of level (0,0) pants}).

\subsection{The level $(0,0)$ tube and cap}\label{subsec: level 0 tube and
cap} We complete the proof of Theorem \ref{thm: gluing defines TQFT
functor} by calculating the series 
$\mathsf{GW}(0\vertline  0,0)_\alpha^\beta$.

\begin{lemma}\label{lem: charge (0,0) tube}
The invariants of the level $(0,0)$ tube are given by:
\[
\mathsf{GW}(0\vertline  0,0)_{\alpha \beta }=\begin{cases}
\frac{1}{\combinatfactor (\alpha ) 
(t_{1}t_{2})^{l (\alpha )}} &\text{if $\alpha =\beta $}\\
0&\text{if $\alpha \neq \beta $.}
\end{cases}
\]
Consequently, we have 
\[
\mathsf{GW}(0\vertline  0,0)_{\alpha }^{\beta }=\delta _{\alpha }^{\beta }
\]
as was required for $\mathsf{GW}(-)$ to be a functor.
\end{lemma}
\textsc{Proof:} 
The virtual dimension of the moduli space
$\overline{M}_{h} (\P ^{1},\alpha ,\beta )$ with {\em connected} domains is
\[
2h-2+l (\alpha )+l (\beta ).
\]
Let $\E^\vee$ be the rank $h$ dual Hodge bundle on
$\overline{M}_h(\P^1,\alpha,\beta)$.  Since the line bundles $L_i$ may be
taken to be trivial,
$$c(-R^\bullet \pi_*f^*(L_i))=c(\E^\vee),$$
where the equality is of ordinary (non-equivariant) Chern classes.
The integral
\begin{multline}
\label{ppw3}
\int _{[\overline{M}_{h} (\P ^{1},\alpha ,\beta )]^{vir}}c_{b_{1}}
(  -R^\bullet \pi_*f^*(L_1)      )c_{b_{2}} (  -R^\bullet \pi_*f^*(L_2)   )= \\
\int _{[\overline{M}_{h} (\P ^{1},\alpha ,\beta )]^{vir}}c_{b_{1}}
(\E^{\vee })c_{b_{2}} (\E^{\vee })
\end{multline}
is zero if 
\[
2h-2+l (\alpha )+l (\beta )>2h.
\]
The only possible non-zero integrals are for $l (\alpha )=l (\beta )=1$.
For $h>0$, 
$$c_{h}(\E^\vee)^2=0,$$
by Mumford's relation. Hence, the
integral \eqref{ppw3} is zero unless $h=0$.

Therefore, the only connected stable map which contributes to the
integral \eqref{ppw3} is the unique degree $d$ map
\[
f_{d}:\P ^{1}\to \P ^{1}
\]
totally ramified over $0$ and $\infty $. The only
disconnected maps which contribute are disjoint unions of genus 0
totally ramified maps of lower degree.
Given a partition $\alpha \vdash d$, let 
\[
f_{\alpha }:\bigsqcup _{l (\alpha )}\P ^{1}\to \P ^{1}
\]
be the map determined by $f_{\alpha _{i}}$ on the  $i$th component. The map
$f_{\alpha }$ has ramification profile $\alpha $ over both $0$ and $\infty$. 
The map is isolated in moduli and has an automorphism group of order
$\combinatfactor (\alpha )$. Thus
\[
\mathsf{GW}^{b_{1},b_{2}} (0\vertline  0,0)_{\alpha\beta}
=\begin{cases}
\frac{1}{\combinatfactor (\alpha )}&\text{if $b_{1}=b_{2}=0$ and $\alpha =\beta $}\\
0&\text{otherwise.}
\end{cases}
\]
The Lemma then follows directly from equation~(\ref{eqn: expression for
Zdgk1k2_lambdas}).\qed

The level $(0,0)$ cap has a simple form obtained by a similar dimensional
argument.

\begin{lemma}\label{lem: charge (0,0) cap}
The invariants of the level $(0,0)$ cap are given by
\[
\mathsf{GW} (0\vertline  0,0)_{\lambda }=\begin{cases}
\frac{1}{d!(t_{1}t_{2})^{d}}& \text{if $\lambda= (1^{d})$}\\
0&\text{if $\lambda \neq (1^{d})$.}
\end{cases}
\]
\end{lemma}
\textsc{Proof:} The (connected domain) moduli space $\overline{M}_{h} (\P
^{1},\lambda )$ has virtual dimension
\[
2h-2+d+l (\lambda ).
\]
Hence,
\begin{equation}\label{eqn: charge 0 cap, connected integral}
\int _{[\overline{M}_{h} (\P ^{1},\lambda )]^{vir}}c_{b_{1}}
(\E^{\vee })c_{b_{2}} (\E^{\vee }) = 0
\end{equation}
if 
\[
2h-2+d+l (\lambda )>2h.
\]
In order for (\ref{eqn: charge 0 cap, connected integral}) to be non-zero, we
must have $d=l (\lambda )=1$. The virtual dimension is then $2h$, which
implies $h=0$ by Mumford's relation.

The only connected stable map for which the integral
(\ref{eqn: charge 0 cap, connected integral}) is non-zero is the
isomorphism 
\[
f:\P ^{1}\to \P ^{1}.
\]
The Lemma is then obtained from \eqref{eqn: expression for
Zdgk1k2_lambdas} by accounting for disconnected
covers.\qed

\subsection{The Calabi-Yau cap}
\begin{lemma}\label{lem: the CY cap}
The invariants of the level $(-1,0)$ cap are given by
\[
\mathsf{GW} (0\vertline  -1,0)_{\lambda }= (-1)^{|\lambda|} 
(-t_{2})^{-l (\lambda
)}\frac{1}{\combinatfactor (\lambda )}\prod _{i=1}^{l (\lambda
)}\left(2\sin \frac{\lambda _{i}u}{2} \right)^{-1}
\]
\end{lemma}
\textsc{Proof:} The calculation has already been done by
localization in the proof of Theorem 5.1 in \cite{Br-Pa-TQFT} in case
$t_1=t_2=1$. We must insert the equivariant parameters. The
relevant {\em connected} integrals are
\[
\int _{[\overline{M}_{h} (\P ^{1},\lambda )]^{vir}}c_{b_{1}} (-R^{\bullet
}\pi _{*}f^{*}\O (-1))c_{b_{2}} (-R^{\bullet }\pi _{*}f^{*}\O ).
\]
The virtual dimension of the moduli space $\overline{M}_{h} (\P
^{1},\lambda )$ is
\[
2h-2+d+l (\lambda ).
\]
The object $-R^{\bullet }\pi _{*}f^{*}\O (-1)$ is represented by a
bundle of rank $h-1+d$. Similarly, $-R ^{\bullet }\pi _{*}f^{*}\O $ is
represented by a bundle of rank $h$ (minus a trivial
factor). Consequently, the integral is zero unless $b_{1}=h-1+d$,
$b_{2}=h$, and $\lambda = (d)$. From equation~(\ref{eqn: expression
for Zdgk1k2_lambdas}), we find that the insertion of the equivariant
parameters yields a factor of $t_{2}^{-1}$.

Since the disconnected invariant
is a product of $l (\lambda )$ connected integrals, the invariant has the
factor $t_{2}^{-l (\lambda )}$.
\qed
\vspace{+10pt}

The series $\mathsf{GW}(0\vertline  0,-1)_\lambda$ 
is obtained from $\mathsf{GW}(0\vertline  -1,0)_\lambda$
by exchanging $t_1$ and $t_2$.

\subsection{The level $(0,0)$ pair of pants}

\subsubsection{Normalization}\label{n3451}
We will study the local theory of curves here with a slightly different
normalization. Recall 
\[
\delta =\sum _{i=1}^{r} (d-l (\lambda ^{i})).
\]
Let
\begin{align}\label{eqn: Zhat}
{\mathsf{GW}}^*
(g\vertline  k_{1},k_{2})_{\lambda ^{1}\dots \lambda ^{r}}&= 
(-iu)^{d (2-2g+k_{1}+k_{2})-\delta } \ 
{\mathsf Z}'(N)_{ \lambda^1 \dots \lambda^r}\\
&= \ \  (-i)^{d (2-2g+k_{1}+k_{2})-\delta }\ \mathsf{GW} 
(g\vertline  k_{1},k_{2})_{\lambda ^{1}\dots \lambda ^{r}}.\nonumber 
\end{align}
With the altered metric,
\[
{\mathsf{GW}}^* (g\vertline k_{1},k_{2})_{\mu ^{1}\dots \mu ^{s}}^{\nu
^{1}\dots \nu ^{t}}= \left(\prod _{i=1}^{t} \combinatfactor (\nu ^{i})
(-t_{1}t_{2})^{l (\nu ^{i})} \right) {\mathsf{GW}}^* (g\vertline
k_{1},k_{2})_{\mu
^{1}\dots \mu ^{s}\nu ^{1}\dots \nu ^{t}},
\]
the partition functions \eqref{eqn: Zhat} satisfy the same gluing rules
as partition functions $\mathsf{GW}
(g\vertline  k_1,k_2)_{\lambda^1\dots \lambda^r}$.
Moreover, a tensor functor,
\[
{\mathbf{GW^*}}:2\mathbf{Cob}^{L_{1},L_{2}}\to R\mathbf{mod}.
\]
is defined just as before.

The reason for the altered normalization is the following
result proven in the Appendix.

\begin{theorem}\label{thm: Z is a rational function of q}
The product
\[
e^{\frac{idu}{2} (2-2g+{k_{1}}+{k_{2}})}\  
\mathsf{GW}^*(g\vertline  k_{1},k_{2})_{\lambda ^{1}\dots \lambda ^{r}} ,
\]
is a \emph{rational function} of $t_{1}$, 
$t_{2}$, and  $q=-e^{iu}$ with $\qnums$-coefficients.
\end{theorem}

The Theorem is closely related to the GW/DT correspondence discussed in
Section~\ref{subsec: GW/DT for local curves}.
The Calabi-Yau cap provides a good example:
\begin{eqnarray*}
e^{\frac{idu}{2}}\  
\mathsf{GW}^*
(0\vertline  -1,0)_{\lambda}& =& e^{\frac{idu}{2}} (-i)^{l(\lambda)}
  (-1)^{d} (-t_{2})^{-l (\lambda
)}\frac{1}{\combinatfactor (\lambda )}\prod _{j=1}^{l (\lambda
)}\left(2\sin \frac{\lambda _{j}u}{2} \right)^{-1} \\
& = & 
(-1)^{d-l(\lambda)} \frac{1}{\combinatfactor (\lambda )}\frac{1}
{ t_2^{l(\lambda)}}   
\prod _{j=1}^{l (\lambda
)} \frac{1}{1 - (-q)^{-\lambda_j}}.
\end{eqnarray*}

\subsubsection{The degree 1 case}
The level $(0,0)$ tube and cap in degree 1 are:
$$\mathsf{GW}^*
(0\vertline  0,0)_{\tableau{1},\tableau{1}} = -\frac{1}{t_1t_2}, 
\ \mathsf{GW}^*(0\vertline  0,0)_{\tableau{1}} = -\frac{1}{t_1t_2}.$$
By the gluing formula,
$$\mathsf{GW}^*
(0\vertline  0,0)_{\tableau{1},\tableau{1},\tableau{1}}
\ (-t_1t_2)\ \mathsf{GW}^*(0\vertline  0,0)_{\tableau{1}} =
\mathsf{GW}^*(0\vertline  0,0)_{\tableau{1},\tableau{1}}.$$
We conclude,
$$\mathsf{GW}^*
(0\vertline  0,0)_{\tableau{1},\tableau{1},\tableau{1}} = -\frac{1}{t_1t_2}.$$
Hence, all the basic series in degree 1 are known.

\subsubsection{The series 
$\mathsf{GW}^*(0\vertline  0,0)_{(d),(d),(1^{d-2}2)}$}
\label{qw223}
The degree $d\geq 2$ series 
$\mathsf{GW}(0\vertline  0,0)_{(d),(d),(1^{d-2}2)}$ plays a
special role in the level $(0,0)$ theory. 

\begin{theorem}\label{thm: the dd2 integral}
For $d\geq 2$,
\begin{eqnarray*}
\mathsf{GW}^* (0\vertline  0,0)_{(d), (d),
(1^{d-2}2)}& = &-\frac{i}{2}\frac{t_{1}+t_{2}}{t_{1}t_{2}}\left(d\cot
\left(\frac{du}{2}\right)-\cot \left(\frac{u}{2}\right) \right) \, .
\end{eqnarray*}
\end{theorem}

\textsc{Proof:} We abbreviate the partition $(1^{d-2}2)$ by $(2)$. After
adjusting
equation \eqref{eqn: expression for Zdgk1k2_lambdas} for the new normalization, we
find
\begin{multline*}
\mathsf{GW}^* 
(0\vertline  0,0)_{(d), (d), (2)}=\\ 
-\frac{i} {(t_{1}t_{2})^{\frac{1}{2}}} \sum
_{b_{1},b_{2}=0}^{\infty }u^{b_1+b_2}\left(\frac{t_{1}}{t_{2}}
\right)^{\frac{b_{2}-b_{1}}{2}} \int _{[\M (\P ^{1},(d), (d), (2))]^{vir}}
c_{b_{1}} (\E^{\vee })c_{b_{2}} (\E^{\vee }).
\end{multline*}
The domains of the maps in the moduli space $\M (\P ^{1},(d), (d),
(2))$ are necessarily connected since there exists a point of total
ramification. Since the virtual dimension is
\[
\operatorname{virdim}\overline{M}_{h} (\P ^{1},(d), (d), (2))=2h-1,
\]
the only values of $(b_{1},b_{2})$ which contribute to 
$\mathsf{GW}^*(0\vertline  0,0)_{(d), (d), (2)}$ are 
$(h,h-1)$ and $(h-1,h)$. We obtain:
\[
\mathsf{GW}^* 
(0\vertline  0,0)_{(d), (d), (2)}=-i\frac{t_{1}+t_{2}}{t_{1}t_{2}}\sum
_{h=1}^{\infty }u^{2h-1}\int _{[\overline{M}_{h} (\P ^{1},(d), (d),
(2))]^{vir}}\rho ^{*} (-\lambda _{h}\lambda _{h-1}).
\]
Here, $\lambda _{k}$ is the $k^{th}$ Chern class of the Hodge bundle on
$\overline{M}_{h,2}$, and
\[
\rho :\overline{M}_{h} (\P ^{1},(d), (d), (2))\to \overline{M}_{h,2}
\]
is the natural map which takes a relative stable map to the domain marked
by the two totally ramified points.

Let 
$
H_{d}\subset M_{h,2}
$
be the locus of curves admitting a degree $d$ map to $\P ^{1}$ which is totally
ramified at the marked points. Equivalently, $H_{d}$ is the locus of curves
$$[C,x_{1},x_{2}]$$ for which $\O (x_{1}-x_{2})$ is a nonzero $d$-torsion point in
$\operatorname{Pic}^{0} (C)$. Let $$\overline{H}_{d}\subset
\overline{M}_{h,2}$$ be the closure of $H_{d}$.

Consider the locus of maps with nonsingular domains,
 $$M_h(\P^1,(d),(d),(2))\subset \overline{M}_h(\P^1,(d),(d),(2)),$$
and let 
$$\partial \overline{M}_h(\P^1,(d),(d),(2))$$
denote the complement.
Let
$$\partial \overline{M}_{h,1} \subset \overline{M}_{h,1}$$
denote the nodal locus.
Let $$\epsilon: \overline{M}_{h,2} \rightarrow\overline{M}_{h,1}$$
be the map forgetting the first point.
An elementary argument yields
$$\rho\left( \partial \overline{M}_h(\P^1,(d),(d),(2)) \right) \subset \epsilon^{-1}(\partial
\overline{M}_{h,1}).$$

The restriction of the virtual class to
$M_h(\P^1,(d),(d),(2))$ is well-known to equal the ordinary fundamental
class of the moduli space, see \cite{pandharipande-toda2000}. Since
$$\rho:M_h(\P^1,(d),(d),(2)) \rightarrow H_d$$
is a proper cover of degree $2h$,
we conclude
\begin{equation}\label{eqn: pushforward of dd2 rel space to Mh2}
\rho _{*}[\overline{M}_{h} (\P ^{1},(d), (d) ,(2))]^{vir}=2h[\overline{H}_{d}]+B
\end{equation}
where $B$ is a cycle supported on $\epsilon^{-1}(\partial
\overline{M}_{h,1})$.

Since $\lambda _{h}\lambda _{h-1}$ vanishes on cycles supported on the
boundary of $\overline{M}_{h,1}$, we find
\[
\mathsf{GW}^* 
(0\vertline  0,0)_{(d), (d), (2)}=i\frac{t_{1}+t_{2}}{t_{1}t_{2}}\sum
_{h=1}^{\infty }u^{2h-1} c_h(d),
\]
where
\[
c_{h} (d)=2h\int _{[\overline{H}_{d}]}\lambda _{h}\lambda _{h-1}.
\]

The cycle $[H_{d}]$ can be described as follows. Let
\[
\begin{diagram}[height=0.8cm]
\operatorname{\mathcal{P}ic}^{0}\\
\uTo^{s} \dTo_{\pi }\\
M_{h,2}
\end{diagram}
\]
be the universal Picard bundle with section
\[
s:[C,x_{1},x_{2}]\mapsto \O (x_{1}-x_{2}).
\]
Let $P_{d}\subset \operatorname{\mathcal{P}ic^{0}}$ be the locus of
\emph{nonzero} $d$-torsion points. Then, by our previous characterization
of $H_{d}$, 
\[
[H_{d}]=\pi _{*}\left(s_{*}[M_{h,2}]\cap P_{d} \right)\in A_*(M_{h,2}).
\]
By a result of Looijenga using the Fourier-Mukai transform, the locus
of $d$-torsion points of {\em any} family of Abelian varieties is a
multiple of the zero section in the Chow ring
\cite{Looijenga-invent}. Hence,
\[
[P_{d}]=\frac{d^{2h}-1}{2^{2h}-1}[P_{2}]
\]
and 
\[
[H_{d}]=\frac{d^{2h}-1}{2^{2h}-1}[H_{2}].
\]
We conclude
\[ c_{h}
(d)=\frac{d^{2h}-1}{2^{2h}-1}c_{h} (2).
\]

Consider the $d=2$ case.
In genus 1, the class $$[\overline{H}_{2}]\in A_*(\overline{M}_{1,2})$$
pushes forward to $3[\overline{M}_{1,1}]$ under the map $\epsilon$.
For genus $h>1$, let $\overline{H}\subset \overline{M}_{h}$ denote the
hyperelliptic locus. There are 
$$(2h+2) (2h+1)$$ ways of marking two of the
Weierstrass points on each curve in $H$. Consequently, the class
$$[\overline{H}_{2}]\in A_*( \overline{M}_{h,2})$$ pushes forward to $(2h+2)
(2h+1)[\overline{H}]$ under the forgetful map 
$$\overline{M}_{h,2}\to
\overline{M}_{h}.$$
 We find
\begin{align*}
\mathsf{GW}^* (0\vertline  0,0)_{(2), (2), (2)}&
=i\frac{t_{1}+t_{2}}{t_{1}t_{2}}\left(6u\int
_{\overline{M}_{1,1}}\lambda _{1}+\sum _{h=2}^{\infty } \frac{(2h+2)!}{(2h-1)!}
u^{2h-1}\int _{\overline{H}}\lambda _{h}\lambda _{h-1}\right)\\
&=i\frac{t_{1}+t_{2}}{t_{1}t_{2}}\left(\frac{u^{4}}{96}+
\sum _{h=2}^{\infty }u^{2h+2}\int _{\overline{H}}\lambda _{h}\lambda _{h-1} \right)'''\\
&=i\frac{t_{1}+t_{2}}{t_{1}t_{2}}\left(u^{2}H (u) \right)'''
\end{align*}
where $H (u)$ is defined in \cite{Faber-Pandharipande-logarithmic} on page 222. 
By Corollary 2 of
\cite{Faber-Pandharipande-logarithmic}, 
\[
(u^{2}H (u))''=-\log \left(\cos \left(\frac{u}{2} \right) \right),
\]
and thus
\[
\mathsf{GW}^* 
(0\vertline  0,0)_{(2), (2), (2)}=\frac{i}{2}
\frac{t_{1}+t_{2}}{t_{1}t_{2}}\tan
\left(\frac{u}{2} \right).
\]
We conclude 
\[
\sum _{h=1}^{\infty }c_{h} (2)u^{2h-1}=\frac{1}{2}\tan\left(\frac{u}{2}
\right).
\]

The function $\cot\left(\frac{u}{2} \right)$ 
is an odd series in $u$ with a simple
pole at $u=0$. We define $b_{h}$ by
\[
\cot\left(\frac{u}{2} \right) = \sum _{h=0}b_{h}u^{2h-1}.
\]
The identity
\[
\frac{1}{2}\tan\left(\frac{u}{2} \right) =
\frac{1}{2}\cot\left(\frac{u}{2} \right)-\cot\left(2\frac{u}{2} \right)
\]
yields
\[
c_{h} (2)= (\frac{1}{2}-2^{2h-1})b_{h}.
\]
Hence,
\[
c_{h} (d)=\frac{1}{2} (1-d^{2h})b_{h}.
\]
We obtain 
\[
\mathsf{GW}^* (0\vertline  0,0)_{(d), (d),
(2)}=-\frac{i}{2}\frac{t_{1}+t_{2}}{t_{1}t_{2}}\left(
d\cot\left(\frac{du}{2} \right)-     \cot \left(\frac{u}{2}
\right)               \right)
\]
which concludes the proof.\qed 

We may write the series as a rational function in $-q=e^{iu}$,
\begin{equation}\label{fferr}
\mathsf{GW}^* (0\vertline  0,0)_{(d), (d),
(2)} =
 \frac{1}{2} \frac{t_{1}+t_{2}}{t_{1}t_{2}} \left(
d\frac{(-q)^d+1}{(-q)^d-1} - \frac{(-q)+1}{(-q)-1}\right).
\end{equation}

\subsection{Reconstruction for the level $(0,0)$ pair of pants}
The main result proven in the Appendix (with C. Faber and A. Okounkov) is the
following.
\begin{theorem}\label{thm: reconstruction of level (0,0) pants}
Let $d\geq 2$.  The set of degree $d$, level $(0,0)$ pair of pants series
\[
\mathsf{GW}^*(0\vertline  0,0)_{\lambda \mu  \nu}
\]
can be uniquely reconstructed from
\[
\mathsf{GW}^* (0\vertline  0,0)_{(d), (d),
(2)}
\]
via the TQFT associativity relations, lower degree series of level $(0,0)$,
and Hurwitz numbers of covering genus $0$.
\end{theorem}

The proof yields an effective method of computing the level $(0,0)$ pair of
pants series via recursions in degree.  Since all the basic series
\eqref{eqn: basic series} can be computed, the full local theory of curves is
effectively determined.  Theorem \ref{thm: Z is a rational function of q}
is obtained in the Appendix as a Corollary of Theorem~\ref{thm:
reconstruction of level (0,0) pants}.

\section{The anti-diagonal action}\label{sec: s1=-s2 limit}
\subsection{Overview}
We study a well-behaved special case of the local
theory of curves. Consider the action of the 
anti-diagonal subgroup
$$\TT = \{ (\xi, \xi^{-1})\   | \ \xi \in \cnums^*\} \subset T.$$
on $N=L_{1}\oplus L_{2}$. The anti-diagonal action 
corresponds
to the limit
\[
t_{1}+t_2=0
\]
in equivariant cohomology.
The 
induced $\TT$-action on $K_{N}$ is trivial.
Explicit formulas can be found since the level
$(0,0)$ Frobenius algebra can be explicitly diagonalized in
the anti-diagonal case.

We define the $Q$-dimension of $\rho $, an irreducible representation
of the symmetric group, indicated $\dim _{Q}\rho $, as follows:
\[
\frac{\dim _{Q}\rho }{d!}= 
\prod _{\Box \in \rho } i\left({ Q^{\frac{h (\Box)}{2}}
-Q^{\frac{-h (\Box)}{2}}} \right)^{-1},
\]
see \cite{Okounkov-Pandharipande-unknot}.
Under the substitution $Q=e^{iu}$, the $Q$-dimension can be expressed as:
\[
\frac{\dim _{Q}\rho }{d!}
= \prod _{\Box \in \rho } 
\left(2 \sin\frac{ h (\Box)u}{2} \right)^{-1}.
\]
By the hook length formula for $\dim \rho $, the
leading term in $u$ of the above expression is $\frac{\dim \rho}{d!} $.

The main result 
here is a closed formula for the (absolute) local theory of curves 
with the anti-diagonal action.

\begin{theorem}\label{thm: s1=-s2 general formula}
Under the restrictions $t_{1}=t$ and $t_{2}=-t$,
\begin{multline*}
\mathsf{GW}_d(g\vertline  k_{1},k_{2})
=\\ (-1)^{d (g-1-k_{2})}t^{d (2g-2-k_{1}-k_{2})}\
\sum
_{\rho }\left(\frac{d!}{\dim \rho } \right)^{2g-2}\left(\frac{\dim \rho
}{\dim _{Q}\rho } \right)^{k_{1}+k_{2}} Q^{\frac{1}{2}c_{\rho }
(k_{1}-k_{2})}
\end{multline*}
where $Q=e^{iu}$ and $c_{\rho }$ is the total content of $\rho $
(see Section \ref{subsub: conventions for partitions}).
\end{theorem}

\subsection{Corollaries}

If $k_{1}+k_{2}=2g-2$, the threefold $N=L_{1}\oplus L_{2}$ is
Calabi-Yau. As previously remarked, the $t_{1}+t_{2}=0$ 
limit corresponds to
the  trivial $\TT$-action on the canonical
bundle. In other words, $N$ is \emph{equivariantly} Calabi-Yau. 

\begin{cor}\label{cor: s1=-s2 limit, CY case}
In the equivariantly Calabi-Yau case,
\[
\mathsf{GW}_d (g\vertline  k,2g-2-k)
= (-1)^{d (g-1-k)}\sum _{\rho }\left(\frac{d!}{\dim
_{Q}\rho } \right)^{2g-2} Q^{-c_{\rho } (g-1-k)}.
\]
\end{cor}
In particular, for the balanced splitting,
$$k_1=k_2=g-1,$$ 
the partition function is a
$Q$-deformation of the classical formula for unramified covers.
\begin{cor}\label{cor: s1=-s2, generic CY}
In the balanced equivariantly Calabi-Yau case,
\[
\mathsf{GW}_d (g\vertline  g-1,g-1)
=\sum _{\rho }\left(\frac{d!}{\dim _{Q}\rho } \right)^{2g-2}.
\]
\end{cor}
Another special Calabi-Yau case is when the base curve $X$ is  elliptic.
We obtain a formula recently derived by Vafa using string theoretic
methods (page 8 of  \cite{Vafa-04-2dYang-Mills}).
\begin{cor}\label{cor: vafa's formula}
Let $L\to E$ be a degree $k$ line bundle on an elliptic curve $E$. The
partition function for the Calabi-Yau action on $L\oplus L^{-1}$ is 

\[
\mathsf{GW}_d(1\vertline  k,-k)=(-1)^{dk}\sum _{\rho } Q^{kc_{\rho }}.
\]
\end{cor}

\subsection{Proof of Theorem \ref{thm: s1=-s2 general formula}}

To derive the formula of Theorem~\ref{thm: s1=-s2 general formula}, we
first explicitly diagonalize the level $(0,0)$ Frobenius algebra for the
anti-diagonal action.
\begin{lemma}\label{lem: s1=-s2 charge 0,0 structure constants}
For the anti-diagonal action, the level $(0,0)$ series have no
nonzero terms of positive degree in $u$. 
\end{lemma}

\textsc{Proof of Lemma.} Let $\cnums _{t}$ denote  $\TT$-representation
given by the standard action of the projection of
$$\TT \subset T= \cnums^* \times \cnums^*$$
on the first factor, and let
\[
c_{1} (\cnums _{t})=t.
\]
The dual line bundle is $\cnums _{t}^{\vee }=\cnums
_{-t}$. 

The level $(0,0)$ partition functions are built from the
following integrals:
\[
\int _{[\M (X,\lambda ^{1},\dots ,\lambda ^{r})]^{vir}}
e \left(-R^{\bullet }\pi
_{*} (\O \otimes \cnums _{t}) \right)
e (-R^{\bullet }\pi
_{*} (\O \otimes \cnums _{-t})).
\]
For any vector bundle $E$, the equivariant Euler class $e (E\otimes \cnums
_{t})$ is a polynomial in $t$ whose coefficients are the (ordinary) Chern
classes of $E$. The above integrand is a weight factor times
\[
e (\E ^{\vee }\otimes \cnums _{t})e (\E^{\vee } \otimes \cnums _{-t}) =
(-1)^{h} e \left(\left(\E ^{\vee }\oplus \E \right)\otimes \cnums _{t}
\right).
\]
Since the Chern classes of $\E ^{\vee }\oplus \E $ all vanish by Mumford's
relation, the last expression is pure weight. The  only non-zero
integrals occur when
 $$b_{1}=b_{2}=0$$
in equation \eqref{eqn: expression for Zdgk1k2_lambdas}.
Only the constant terms in $u$ are
non-zero. In particular, $\mathsf{GW} 
(0\vertline  0,0)^{\gamma }_{\alpha \beta }$ is given
by the $t_{1}+t_{2}=0$ limit of equation (\ref{eqn: u=0 structure
constants}).\qed
\vspace{+10pt}

The structure
constants for the level $(0,0)$ Frobenius algebra are given by:
\[
\mathsf{GW} 
(0\vertline  0,0)_{\alpha \beta }^{\gamma }=\left(-t^{2} \right)^{\frac{1}{2}
(d-l (\alpha )-l (\beta )+l (\gamma ))}\sum _{\rho }\left(\frac{d!}{\dim
\rho } \right)\frac{\chi^{\rho }_\alpha \chi^{\rho }_\beta}
{\combinatfactor (\alpha )\combinatfactor (\beta )}\chi^{\rho}_
\gamma.
\]
As a consequence of the Lemma, multiplication in the level $(0,0)$
Frobenius algebra is diagonalized by the basis $v_{\rho }^{0}$ constructed
in the proof of Proposition~\ref{prop: Zd is semi-simple}.

In order to diagonalize the level $(0,0)$ Frobenius algebra,
 we had to enlarge the
coefficient ring to $\tilde{R}$ to include the formal square roots
$t_{1}^{\frac{1}{2}}$ and $t_{2}^{\frac{1}{2}}$. The specialization
\[
t_{1}^{\frac{1}{2}}=t^{\frac{1}{2}},\quad t_{2}^{\frac{1}{2}}
=it^{\frac{1}{2}}
\]
is compatible with
$$t_1=t, \quad t_2=-t.$$
By Lemma \ref{lem: s1=-s2 charge 0,0 structure constants}, the idempotent
basis $v_{\rho }^{0}=v_{\rho }$ given by equation \eqref{eqn: v0 basis} is:
\begin{equation}\label{dfg}
v_{\rho }=\frac{\dim \rho }{d!}\sum _{\alpha } (it)^{l (\alpha )-d}\chi
^{\rho }_\alpha e_{\alpha }.
\end{equation}

To apply Theorem \ref{thm: structure theorem}, we must
 compute $\lambda _{\rho }$, $\eta _{\rho }$, and $\overline{\eta }_{\rho
}$. We compute $\lambda _{\rho }$ as follows:
\begin{align*}
\lambda _{\rho }^{-1}&=\mathbf{GW} \left(\intextcounitFA
^{\sss{(0,0)}} \right) (v_{\rho })\\
&=\frac{\dim \rho }{d!}\sum _{\alpha } (it)^{l (\alpha )-d}\chi
^{\rho }_\alpha \ \mathsf{GW} (0\vertline  0,0)_{\alpha }\\
&=\frac{\dim \rho }{d!} (it)^{l (1^{d} )-d}\chi
^{\rho }_{(1^{d})}\frac{1}{d!(-t^{2})^{d}}\\
&=\left(\frac{\dim \rho }{d!} \right)^{2} (it)^{-2d}.
\end{align*}
Hence,
$$\lambda _{\rho } =(it)^{2d}\left(\frac{d!}{\dim \rho } \right)^{2}.$$

In order to compute $\eta _{\rho}$, we must express $\eta $ in 
terms of the basis
$\{v_{\rho } \}$.
\begin{align*}
\eta& =\mathbf{GW}\left(\intextunitFA
^{\sss{(-1,0)}} \right)\\
&=\sum _{\alpha }
\mathsf{GW} (0\vertline  -1,0)^{\alpha }e_{\alpha }\\
&=\sum _{\alpha } (-1)^{d}t^{l (\alpha )}\left(\prod _{i=1}^{l (\alpha )}
\frac{-1}{2\sin \frac{\alpha _{i}u}{2}} \right)e_{\alpha }\\
&=\sum _{\alpha } (-1)^{d} (it)^{l (\alpha )}Q^{d/2}\left(\prod _{i=1}^{l (\alpha )}
\frac{1}{1-Q^{\alpha _{i}}} \right)e_{\alpha }
\end{align*}
where $Q=e^{iu}$ as before. The expression
\[
\prod _{i=1}^{l (\alpha )}\frac{1}{1-Q^{\alpha _{i}}}
\]
arises in the theory of symmetric functions. 
The power sum symmetric functions are defined by:
\begin{align*}
p_{k} (x_{1},x_{2},x_3,\dots )&= x_{1}^{k}+x_{2}^{k}+x_3^k+\dots \\
p_{\alpha }&=\prod _{i=1}^{l (\alpha )}p_{\alpha _{i}}.
\end{align*}
For the specialization 
$$x_{1}=1, \ x_{2}=Q, \  x_{3}=Q^2, \dots,$$
we obtain 
$$p_{k} (Q)= (1-Q^{k})^{-1}.$$ 
Hence,
\[
\eta =\sum _{\alpha } (-1)^{d} (it)^{l (\alpha )}Q^{d/2}p_{\alpha }
(Q)\ e_{\alpha }
\]
and similarly
\[
\overline{\eta } =\sum _{\alpha } (-1)^{d} (it)^{l (\alpha
)}Q^{d/2} (-1)^{l (\alpha )} p_{\alpha } (Q)\ e_{\alpha }.
\]

Inversion of \eqref{dfg} yields the following formula:
\begin{equation}\label{kwqq3}
e_{\alpha }=(it)^{d-l (\alpha )}\sum _{\rho }\frac{d!}{\dim \rho
}\frac{\chi^{\rho }_\alpha}{\combinatfactor (\alpha )}v_{\rho }
\end{equation}
After substituting \eqref{kwqq3} in the expression for $\eta$, we find
\begin{align*}
\eta &=\sum _{\rho }v_{\rho} \left[ (-it)^{d}Q^{d/2}\frac{d!}{\dim
\rho }\left(\sum _{\alpha }\frac{\chi^{\rho }_\alpha p_{\alpha }
(Q)}{\combinatfactor (\alpha )} \right)\right],\\
\overline{\eta } &=\sum _{\rho }v_{\rho} \left[ (+it)^{d}Q^{d/2}\frac{d!}{\dim
\rho }\omega \left(\sum _{\alpha }\frac{\chi^{\rho }_\alpha p_{\alpha }
(Q)}{\combinatfactor (\alpha )} \right)\right].
\end{align*}
Here, $\omega $ is the involution on the ring of symmetric
functions defined by
\[
(-1)^{l (\alpha )}p_{\alpha }= (-1)^{d}\omega (p_{\alpha }).
\]

The sum over $\alpha $ in the
latter expressions for $\eta$ and $\overline{\eta}$
is equal to the Schur function $s_{\rho }
(Q)$, see \cite{Macdonald} page 114. 
We have
\[
\omega (s_{\rho })=s_{\rho '}
\]
where $\rho '$ is the dual representation (or conjugate partition),
\cite{Macdonald} page 42.
Thus, we
obtain
\begin{align*}
\eta _{\rho }&= (-it)^{d}Q^{d/2}\frac{d!}{\dim \rho }s_{\rho }
(Q)\\
\overline{\eta } _{\rho }&= (+it)^{d}Q^{d/2}\frac{d!}{\dim \rho }s_{\rho' }(Q).
\end{align*}
The Schur functions are easily expressed in terms of the
$Q$-dimension. From \cite{Macdonald} page 45, 
\begin{align*}
s_{\rho }&=Q^{n (\rho )}\prod _{\Box \in \rho
}\frac{1}{1-Q^{h (\Box )}}\\
&=Q^{n (\rho )-\frac{1}{2} (n (\rho )+n (\rho ')+d)} (-1)^{d}\prod
_{\Box \in \rho }\left(Q^{h (\Box )/2}-Q^{-h (\Box )/2}
\right)^{-1}\\
&=Q^{-\frac{1}{2} (d+c_{\rho })}i^{d}\ \frac{\dim _{Q}\rho }{d!}.
\end{align*}
We have used \eqref{eqn: identities for total content and hook length} in the
above formulas. We conclude
\begin{align*}
\eta _{\rho }&= (+t)^{d}\; Q^{-\frac{c_{\rho }}{2}}\; \frac{\dim _{Q}\rho }{\dim \rho },\\
\overline{\eta } _{\rho }&= (-t)^{d}\; Q^{+\frac{c_{\rho }}{2}}\; \frac{\dim _{Q}\rho }{\dim \rho }.
\end{align*}
Theorem \ref{thm: s1=-s2 general formula} then follows directly from 
Theorem \ref{thm: structure theorem}.\qed

\section{A degree 2 calculation}\label{sec: explicit formulas}

The partition function $\mathsf{GW} (g\vertline  k_{1},k_{2})$ in 
degree 2 is calculated here. The result was announced previously in
\cite{Br-Pa-TQFT}.

We abbreviate the level $(0,0)$ pair of pants by
\[
\mathsf{GW} (0\vertline  0,0)_{\lambda \mu \nu }=
\mathsf{P}_{\lambda \mu \nu }
\]
and the Calabi-Yau cap by
\[
\mathsf{GW} (0\vertline  -1,0)_{\lambda }=\mathsf{C}_{\lambda }.
\]
We apply the usual convention \eqref{eqn: index
raising formula} for raising indices to 
$\mathsf{P}_{\lambda \mu \nu }$ and $\mathsf{C}_{\lambda }$.  

From the proof of Theorem~\ref{thm: structure theorem}, the partition
function is
\[
\mathsf{GW} 
(g\vertline  k_{1},k_{2})=
\Tr \left(G^{g-1}A^{-k_{1}}\overline{A}^{-k_{2}} \right).
\]
The genus adding operator $G$ and the right annihilation operator $A$ can
be computed in terms of 
$\mathsf{P}_{\lambda \mu \nu }$ and $\mathsf{C}_{\lambda }$ by the
gluing formula.
\begin{align}\label{eqn: A and G in terms of Z and C}
G^{\mu }_{\nu }&=\sum _{\lambda ,\epsilon \vdash d}
\mathsf{P}^{\mu }_{\lambda \epsilon }\mathsf{P}^{\lambda \epsilon }_{\nu },\\
A^{\mu }_{\nu }&=\sum _{\lambda \vdash d}
\mathsf{C}^{\lambda }\mathsf{P}^{\mu }_{\lambda \nu }.\nonumber
\end{align}
We obtain $\overline{A}$ from $A$ by switching $t_{1}$ and $t_{2}$.

$\mathsf{P}_{\lambda \mu \nu }$ is determined recursively by Theorem~\ref{thm:
reconstruction of level (0,0) pants} and 
$\mathsf{C}_{\lambda }$ is given explicitly
by Lemma~\ref{lem: the CY cap}. We list their values for $d\leq 2$:

\begin{equation*}
\begin{array}{lll}
\mathsf{C}_{\tableau{1}}=\frac{1}{t_{2}}\frac{1}{2\sin \frac{u}{2}},
&\quad   \mathsf{C}_{\tableau{1 1}}=
\frac{1}{t_{2}^{2}}\frac{1}{2\left(2\sin \frac{u}{2} \right)^{2}},  &\quad  
\mathsf{C}_{\tableau{2}}=-\frac{1}{t_{2}}\frac{1}{4\sin u},\\ &&\\
\mathsf{P}_{\tableau{1}\tableau{1}\tableau{1}}= \frac{1}{t_{1}t_{2}},
&\quad      
\mathsf{P}_{\tableau{1 1}\tableau{1 1}\tableau{1 1}}=
\frac{1}{2} \frac{1}{(t_{1}t_{2})^{2}},&\quad   
\mathsf{P}_{\tableau{1 1}\tableau{1 1}\tableau{2}}=0,\\ &&\\    
\mathsf{P}_{\tableau{1 1}\tableau{2}\tableau{2}}=
\frac{1}{2} \frac{1}{t_{1}t_{2}}.&\quad
\mathsf{P}_{\tableau{2}\tableau{2}\tableau{2}}=
-\frac{1}{2} \frac{t_1+t_2}{t_{1}t_{2}}\tan\frac{u}{2}.&\quad
\\
\end{array}
\end{equation*}

\vspace{10pt}

For $d=1$, we have 
\begin{align*}
G^{\tableau{1}}_{\tableau{1}}&=
\mathsf{P}^{\tableau{1}}_{\tableau{1}\tableau{1}}
\mathsf{P}^{\tableau{1}\tableau{1}}_{\tableau{1}}=t_{1}t_{2}\\ 
A^{\tableau{1}}_{\tableau{1}}&=
\mathsf{C}_{\tableau{1}}
\mathsf{P}_{\tableau{1}}^{\tableau{1}\tableau{1}}
=t_{1} \left(2\sin \frac{u}{2} \right)^{-1}.
\end{align*}
Hence,
\[
\mathsf{GW}_{1} (g\vertline  k_{1},k_{2})
= (t_{1}t_{2})^{g-1}t_{1}^{-k_{1}}t_{2}^{-k_{2}}\left(2\sin \frac{u}{2} \right)^{k_{1}+k_{2}}.
\]

\vspace{10pt}

For $d=2$, we compute the entries of $G$ and $A$ via \eqref{eqn: A
and G in terms of Z and C} to obtain:

\begin{align*}
G&=
\left(
\begin{smallmatrix} 
4 (t_{1}t_{2})^{2}&\quad \quad &-2 (t_{1}t_{2})^{2} 
(t_{1}+t_{2})\tan\frac{u}{2}\\&\quad \quad &\\
-2 (t_{1}t_{2})(t_{1}+t_{2})\tan\frac{u}{2} & 
&4 (t_{1}t_{2})^{2}+2 (t_{1}t_{2}) (t_{1}+t_{2})^{2}\tan^{2}\frac{u}{2}
\end{smallmatrix} 
\right),\\
&\\
A&=\left(
\begin{smallmatrix}
t_{1}^{2}\left(2\sin \frac{u}{2} \right)^{-2}&&-t_{1}^{2}
t_{2}\left(2\sin u \right)^{-1}\\ &&\\
-t_{1}\left(2\sin u \right)^{-1}
&&t_{1} (t_{1}+t_{2})\left(2\cos \frac{u}{2} \right)^{-2}+
t_{1}^{2}\left(2 \sin \frac{u}{2} \right)^{-2}
\end{smallmatrix} 
\right).
\end{align*}
The matrices $G$, $A$, and $\overline{A}$ mutually commute and so we can
simultaneously diagonalize them to obtain:
\begin{equation}\label{eqn: d=2 full formula}
\mathsf{GW}_{2} (g\vertline  k_{1},k_{2})=
\lambda _{+}^{g-1}\eta _{+}^{-k_{1}}\overline{\eta }_{+}^{-k_{2}}+\lambda _{-}^{g-1}\eta _{-}^{-k_{1}}\overline{\eta }_{-}^{-k_{2}},
\end{equation}
where
\begin{align*}
\lambda _{\pm}&=\frac{t_{1}t_{2}}{(1-q)^{2}}
\left(-\Theta \pm (1+q) (t_{1}+t_{2})\sqrt{\Theta } \right)\\
\eta _{\pm }&=\frac{qt_{1}}{2 (1-q^{2})^{2}}\left((t_{1}-t_{2}) 
(1+q)^{2}-8t_{1}q\pm (1+q)\sqrt{\Theta } \right)\\
\overline{\eta } _{\pm }&=\frac{qt_{2}}{2 (1-q^{2})^{2}}
\left((t_{2}-t_{1}) (1+q)^{2}-8t_{2}q\pm (1+q)\sqrt{\Theta } \right)\\
\Theta &= (t_{1}-t_{2})^{2} (1+q)^{2}+16qt_{1}t_{2}\\
q&=-e^{iu}.
\end{align*}
For the specialization 
\[
t_{1}=t_{2}=t,
\]
the above equations simplify to 
\begin{align*}
\lambda _{\pm }&=\frac{4t^{4}}{1\mp \sin \frac{u}{2}}\\
\eta _{\pm }=\overline{\eta }_{\pm }&=\frac{t^{2}}
{4\sin ^{2}\frac{u}{2}\left(1\mp \sin \frac{u}{2} \right)}.
\end{align*}
Hence,
\begin{align*}
\mathsf{GW}_{2} (g\vertline  k_{1},k_{2})\Big|_{t_{1}=t_{2}=t}
=&t^{2 (2g-2-k_{1}-k_{2})}4^{g-1}
\left(2\sin \frac{u}{2} \right)^{2 (k_{1}+k_{2})} \\
&\cdot \left\{\left(1+\sin \frac{u}{2} \right)^{k_{1}+k_{2}+1-g}
+\left(1-\sin \frac{u}{2} \right)^{k_{1}+k_{2}+1-g} \right\}.
\end{align*}
In particular, the local $d=2$ Calabi-Yau partition function 
is given by:
\begin{multline*}
\mathsf{GW}_{2} (g\vertline g-1,g-1)\Big|_{t_1=t_2=t}=\\
\left(2\sin \frac{u}{2} \right)^{4g-4}\left\{\left(4-4\sin
\frac{u}{2} \right)^{g-1}+ \left(4+4\sin \frac{u}{2} \right)^{g-1} \right\},
\end{multline*}
in agreement with the announcement of \cite{Br-Pa-TQFT} up to
$u$ shifting conventions.

The partition function in degree 2
is easily seen to satisfy the BPS integrality of
Gopakumar-Vafa.

\section{The GW/DT correspondence for residues}\label{sec: GW/DT
correspondence for residues}

\subsection{Overview} A
Gromov-Witten/Donaldson-Thomas correspondence parallel to
\cite{MNOP1,MNOP2} is conjectured here for equivariant residues in both the
absolute and relative cases. The computation of the full local
Gromov-Witten theory of curves together with the GW/DT correspondence
predicts the full local Donaldson-Thomas theory of curves.

\subsection{Residue Invariants in Donaldson-Thomas theory}

Let $Y$ be a nonsingular, {\em quasi-projective}, algebraic threefold.
Let $I_{n} (Y,\beta )$ denote the moduli space of ideal sheaves
\[
0 \rightarrow I_{Z}\rightarrow  \O _{Y}\rightarrow \O_Z \rightarrow 0
\]
of subschemes $Z$ of degree $\beta =[Z]\in H_{2}
(Y,\znums )$ and Euler characteristic $n=\chi (\O _{Z})$. 
Though $Y$ may not be compact, we require  $Z$ to have \emph{proper} support.

Let $Y$ be equipped with an action by an algebraic torus $T$.
The moduli space $I_n(Y,\beta)$ carries a $T$-equivariant perfect obstruction 
theory obtained from (traceless) $\text{Ext}_0(I,I)$, see
\cite{t}. Though $Y$ is
quasi-projective, 
$\text{Ext}_0(I,I)$ is well-behaved since the associated
quotient scheme $Z\subset Y$ is proper. Alternatively, for any
$T$-equivariant compactification,
$$Y\subset \overline{Y},$$
the obstruction theory on 
$$I_n(Y,\beta)\subset I_n(\overline{Y},\beta)$$
is obtained by 
restriction.

We will define
Donaldson-Thomas residue invariants under the following assumption.
\begin{assumption}\label{assuption: T fixed locus of DT space is compact}
The $T$-fixed point set $I_{n} (Y,\beta )^{T}$ is compact.
\end{assumption}

The definition of the Donaldson-Thomas residue invariants of $Y$ follows
the 
strategy of the Gromov-Witten case. We define ${\mathsf Z}_{DT}(Y)_{\beta
} $ formally by:
\begin{equation}\label{eqn: moral defn of Zdt}
{\mathsf Z}_{DT}(Y)_{\beta }\text{ ``}=\text{'' }\sum _{n\in \znums
}q^{n} \int _{[I_{n} (Y,\beta ) ]^{vir}}1 .
\end{equation}
The variable $q$  indexes the Euler number $n$.
Under Assumption~\ref{assuption: T fixed locus of DT space is compact},
the integral on
the right of \eqref{eqn: moral defn of Zdt} is well-defined
by the virtual localization formula as an equivariant residue.

\begin{definition} The partition function for
the degree $\beta $ Donaldson-Thomas residue invariants of $Y$ is defined by:
\begin{equation}
\label{defn: definition of Zdt}
{\mathsf Z}_{DT}(Y)_{\beta }  =\sum _{n\in \znums }q^{n}
\int _{[I _{n} (Y,\beta )^{T}]^{vir}}\frac{1}{e
(\operatorname{Norm}^{vir})}.
\end{equation}
\end{definition}

The $T$-fixed
part of the perfect obstruction theory for $I_n(Y,\beta )$ induces a
perfect obstruction theory for $I_n (Y,\beta )^{T}$ and hence a virtual
class \cite{Gr-Pa,MNOP1}. 
The 
equivariant virtual normal bundle
of the embedding, $$I_n (Y,\beta )^{T}\subset I_n (Y,\beta ), $$ is
 $\operatorname{Norm}^{vir}$  with
equivariant Euler class
$e(\operatorname{Norm}^{vir})$. The integral in \eqref{defn: definition of Zdt}
denotes equivariant push-forward to a point.

As defined, ${\mathsf Z}_{DT}(Y)_{\beta }$ is {\em unprimed} since
 the degree $0$ contributions have not yet been removed. In
Gromov-Witten theory, the degree 0 contributions are removed geometrically by
forbidding such components in the moduli problem. Since a geometrical
method of removing the degree 0 contribution from 
Donaldson-Thomas theory does not appear to be available, a formal
method is followed.

\begin{definition}\label{defn: reduced DT partition function Zdt'}
The reduced partition function ${\mathsf Z}'_{DT}(Y)_{\beta }$ for the degree $\beta
$ Donaldson-Thomas residue invariants of $Y$ is defined by:
\begin{equation*}
{\mathsf Z}'_{DT}(Y)_{\beta }  = \frac{{\mathsf Z}_{DT}(Y)_{\beta } }
{{\mathsf Z}_{DT}(Y)_0}. 
\end{equation*}
\end{definition}

Let $r$ be the rank of $T$, and let $t_1,\dots, t_r$ be generators of the
equivariant cohomology of $T$.  By definition, $ {\mathsf
Z}'_{DT}(Y)_{\beta } $ is a Laurent series in $q$ with coefficients given
by rational functions of $t_1,\ldots, t_r$ of homogeneous degree equal to
minus the virtual dimension of $I_n(Y,\beta)$.

\subsection{Conjectures for the absolute theory}

The equivariant degree 0 series is conjecturally determined
in terms of the MacMahon function,
$$M(q)=\prod_{n\geq 1} \frac{1}{(1-q^n)^n},$$
the generating series for 3-dimensional partitions.

\vspace{+10pt}
\noindent {\bf Conjecture 1.} The degree 0 Donaldson-Thomas partition function
is determined by:
\begin{equation*}
{\mathsf Z}_{DT} (Y)_0=M (-q)^{\int _{Y}c_{3} (T_Y\otimes K_Y)},
\end{equation*}
where the integral in the exponent is defined via localization on $Y$,
\[
\int _{Y}c_{3} (T_Y\otimes K_Y) =\int _{Y^{T}}\frac{c_{3}(T_Y\otimes K_Y) }{e (N_{Y^{T}/Y})} \in
{\mathbb Q}(t_1,\dots, t_r).
\]
\vspace{+10pt}

The subvariety $Y^{T} $ is compact as a consequence of
Assumption~\ref{assuption: T fixed locus of DT space is compact}. By
Theorem 1 of \cite{MNOP2}, Conjecture 1 holds for toric $Y$.

\vspace{+10pt}
\noindent {\bf Conjecture 2.} The reduced series ${\mathsf Z}'_{DT}(Y)_{\beta }$
is a rational function of the equivariant parameters $t_i$ and $q$.
\vspace{+10pt}

The GW/DT correspondence for absolute residue invariants can now be stated.

\vspace{+10pt}
\noindent {\bf Conjecture 3.} 
After the change of variables $e^{iu}=-q$, 
\[
(-iu)^{\int_\beta c_1(T_Y)} \ {\mathsf Z}'_{GW}(Y)_{\beta } =(-q)^{-\frac{1}{2}\int_\beta c_1(T_Y)}
\ {\mathsf Z}'_{DT}(Y)_{\beta}.
\]
\vspace{+10pt}

Conjectures 1-3 are equivariant versions of the conjectures of \cite{MNOP1,MNOP2}.
In \cite{MNOP2}, a GW/DT correspondence for primary and (certain) descendent
field insertions is presented. The equivariant correspondence with insertions
remains to be studied.

\subsection{The relative conjectures}
A Gromov-Witten/Donaldson-Thomas residue correspondence for relative theories
may also be defined.
Let 
$$S\subset Y$$
be a nonsingular, $T$-invariant, divisor.
Let
$\beta\in H_2(Y,{\mathbb Z})$
be a curve class
satisfying
$$\int_\beta[S] \geq 0.$$
Let $\eta$ be a partition of $\int_\beta[S]$ weighted
by the equivariant cohomology of $S$, 
$$H_T^*(S,{\mathbb Q}).$$ We  follow here
the notation of \cite{MNOP2}.
The reduced Gromov-Witten partition function,
\begin{equation*}
{\mathsf Z}'_{GW}(Y/S)_{\beta,\eta} =\sum _{h\in \znums }u^{2h-2}\int
_{[\M (Y/S,\beta,\eta )^{T}]^{vir}}\frac{1}{e (\operatorname{Norm}^{vir})},
\end{equation*}
is well-defined if $\M (Y/S,\beta,\eta )^{T}$ is compact.
The weighted partition $\eta$ specifies the relative
conditions imposed on the moduli space of maps.

We refer the reader to \cite{MNOP2} for a discussion of relative
Donaldson-Thomas theory.
The relative Donaldson-Thomas partition function,
\begin{equation*}
{\mathsf Z}_{DT}(Y/S)_{\beta,\eta }  =\sum _{n\in \znums }q^{n}
\int _{[I _{n} (Y/S,\beta,\eta )^{T}]^{vir}}\frac{1}{e
(\operatorname{Norm}^{vir})},
\end{equation*}
is well-defined if $I _{n} (Y/S,\beta,\eta )^{T}$ is compact.
Let
\begin{equation*}
{\mathsf Z}'_{DT}(Y/S)_{\beta,\eta }  = \frac{{\mathsf Z}_{DT}(Y/S)_{\beta,\eta } }
{{\mathsf Z}_{DT}(Y/S)_0}. 
\end{equation*}
denote the reduced relative partition function.

The weighted partition $\eta$ in Donaldson-Thomas theory specifies the relative
conditions imposed on the moduli space of ideal sheaves. The partition
$\eta$
determines an element of the
Nakajima basis of 
the $T$-equivariant cohomology of the Hilbert scheme of
points of $S$.

\vspace{+10pt}
\noindent {\bf Conjecture 1R.} The degree 0 relative Donaldson-Thomas partition function
is determined by:
\begin{equation*}
{\mathsf Z}_{DT} (Y/S)_0=M (-q)^{\int _{Y}c_{3} (T_Y[-S]\otimes K_Y[S])},
\end{equation*}
where $T_Y$ is the sheaf of tangent fields with logarithmic zeros and
$K_Y$ is the logarithmic canonical bundle.
\vspace{+10pt}

\vspace{+10pt} 
\noindent 
{\bf Conjecture 2R.} The reduced series ${\mathsf
Z}'_{DT}(Y/S)_{\beta,\eta }$ is a rational function of the equivariant
parameters $t_i$ and $q$.  
\vspace{+10pt}

\noindent {\bf Conjecture 3R.} 
After the change of variables $e^{iu}=-q$, 
\[
(-iu)^{\int_\beta c_1(T_Y)+l(\eta)-|\eta|} 
\ {\mathsf Z}'_{GW}(Y/S)_{\beta,\eta  } =(-q)^{-\frac{1}{2}\int_\beta c_1(T_Y)}
\ {\mathsf Z}'_{DT}(Y/S)_{\beta,\eta},
\]
where $|\eta|=\int_\beta[S]$.
\vspace{+10pt}

Conjectures 1R-3R are equivariant versions of the relative conjectures of
\cite{MNOP2} without insertions.

\subsection{The local theory of curves}
\label{subsec: GW/DT for local curves}
Let $X$ be a nonsingular curve of genus $g$.
Let $N$ be a rank 2 bundle on a $X$ with a direct sum decomposition,
\[
N=L_{1}\oplus L_{2}.
\]
Let $k_i$ denote the degree of $L_i$ on $X$. 

Consider the Gromov-Witten residue theory of $N$
relative to the
$T$-invariant divisor
$$S= \bigcup_{p\in D} N_{p}\  \subset N,$$
where $D \subset X$ is a finite set of points. Since
$$H_T^*(S) = \bigoplus_{p\in D} H^*_T(p),$$
$\eta$ is simply a list of partitions indexed by $D$.

\begin{theorem} \label{ko34}
The GW/DT correspondence holds (Conjectures 1R-3R)
for the local theory of curves.
\end{theorem}

Theorem \ref{ko34} is proven by matching the calculation of the
local Gromov-Witten theory of curves here with the determination 
of the local Donaldson-Thomas theory of curves in \cite{dtlc}.
The results of \cite{dtlc} depend upon foundational aspects
of relative Donaldson-Thomas theory which have not yet been
treated in the literature.

The $\mathsf{GW}^*$-partition functions in relative Gromov-Witten
theory are defined in Section \ref{n3451}. Similarly, the
$\mathsf{DT}^*$-partition functions in relative Donaldson-Thomas
theory are defined by:
\[
\mathsf{DT}^*(g\vertline k_1,k_2)_{\lambda^1\dots\lambda^r} =
(-q)^{-\frac{d}{2}(2-2g+k_1+k_2)} \ {\mathsf Z}'_{DT}(N)_{d[X],
\lambda^1 \dots \lambda^r}.
\]
The GW/DT correspondence for local curves can be
conveniently restated as the equality,
\begin{equation}\label{kk231}
\mathsf{GW}^*(g\vertline  k_1,k_2)_{\lambda^1\dots\lambda^r} =
\mathsf{DT}^*(g\vertline  k_1,k_2)_{\lambda^1\dots\lambda^r},
\end{equation}
after the variable change $e^{iu}=-q$,

\section{Further directions}\label{sec: further directions}
\subsection{The Hilbert scheme $\Hilb^n(\cnums^2)$}

Let the 2-dimensional torus $T$ act on $\cnums^2$ by scaling the factors.
Consider the induced $T$-action on $\Hilb^n(\cnums^2)$.
The $T$-equivariant cohomology of $\Hilb^n(\cnums^2)$,
$$H^*_T(\Hilb ^{n} (\cnums ^{2}),\qnums),$$
 has
a canonical Nakajima basis, 
$$\{ \ |\mu\rangle\ \}_{|\mu|=n}$$
indexed by partitions of $n$.
The degree of a curve in $\Hilb^n(\cnums^2)$ is determined
by intersection with the divisor
$$D = - |2,1^{n-2} \rangle.$$

Define the series $\langle \lambda,\mu,\nu\rangle^{{\Hilb^n (\cnums^2)}}$ of
 3-pointed, genus 0, $T$-equivariant Gromov-Witten invariants by a sum over curve degrees:
$$\langle \lambda,\mu,\nu
\rangle^{\Hilb^n(\cnums^2)} = \sum_{d\geq 0} q^d
\langle \lambda,\mu,\nu
\rangle^{\Hilb^n(\cnums^2)}_{0,3,d}.$$

\begin{theorem}\label{gg345} A Gromov-Witten/Hilbert correspondence holds:
\begin{equation*} 
\mathsf{GW}^*(0\vertline 0,0)_{\lambda\mu\nu} 
= (-1)^n\langle \lambda,\mu,\nu
\rangle^{\Hilb^n(\cnums^2)},
\end{equation*}
after the variable change $e^{iu}=-q$.
\end{theorem}
The proof of Theorem \ref{gg345} is obtained by our
determination of the series $\mathsf{GW}^*(0\vertline 0,0)_{\lambda\mu\nu}$
together with the computation of the quantum cohomology of
the Hilbert scheme in \cite{Ok-Pan-Hilb}.

\subsection{The orbifold $\Sym (\cnums ^{2}) $} The 3-pointed, genus
0, $T$-equivariant Gromov-Witten invariants of the orbifold $\Sym
(\cnums^2)=(\cnums^2)^n/S_n$ can be related to 
$\mathsf{GW}^*(0\vertline
0,0)_{\lambda\mu\nu}$, see \cite{Bryan-Graber}.

The Hilbert scheme $\Hilb^n(\cnums^2)$ is a crepant resolution of the
(singular) quotient variety $\Sym (\cnums ^{2})$.  Theorem
\ref{gg345} may be viewed as relating the
$T$-equivariant quantum cohomology of the quotient {\em orbifold}
$\Sym (\cnums^2)$ to the $T$-equivariant quantum cohomology of the
resolution $\Hilb^n(\cnums^2)$. The correspondence requires extending
the definition of orbifold quantum cohomology to include quantum
parameters for twisted sectors \cite{Bryan-Graber}.

Mathematical conjectures relating the quantum cohomologies of
orbifolds and their crepant resolutions in the
non-equivariant case have been pursued by 
Ruan (motivated by the physical predictions of Vafa and Zaslow). 
Theorem \ref{gg345} suggests that the correspondence also holds
in the equivariant context.

\vspace{1.5in}
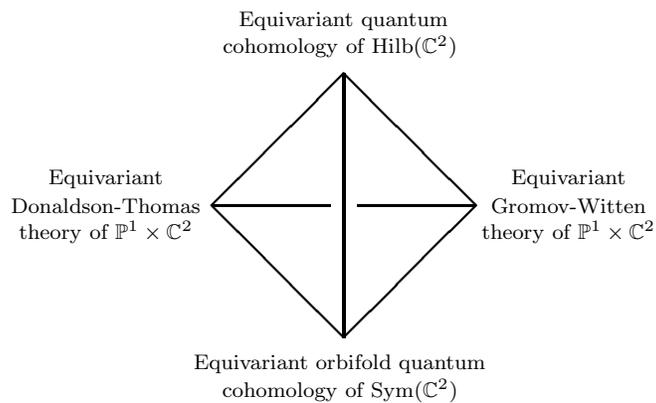
\begin{center}
\scriptsize
\begin{picture}(200,75)(-30,-50)
\put(-100 ,-115 ){\line(1 ,0 ){340}}
\put(-100 ,-115 ){\line(0 ,1 ){225}}
\put(-100 ,110 ){\line(1 ,0 ){340}}
\put(240 ,-115 ){\line(0 ,1 ){225}}
\thicklines
\put(25,25){\line(1,1){50}}
\put(25,25){\line(1,-1){50}}
\put(125,25){\line(-1,1){50}}
\put(125,25){\line(-1,-1){50}}
\put(75,-25){\line(0,1){100}}
\put(25,25){\line(1,0){45}}
\put(80,25){\line(1,0){45}}
\put(75,95){\makebox(0,0){Equivariant quantum}}
\put(75,85){\makebox(0,0){cohomology of $\Hilb (\cnums ^{2})$ }}
\put(75,-35){\makebox(0,0){Equivariant orbifold quantum }}
\put(75,-45){\makebox(0,0){cohomology of $\operatorname{Sym} (\cnums ^{2})$ }}
\put(160,35){\makebox(0,0){Equivariant}}
\put(160,25){\makebox(0,0){Gromov-Witten}}
\put(160,15){\makebox(0,0){theory of $\P ^{1}\times \cnums ^{2}$}}
\put(-15,35){\makebox(0,0){Equivariant}}
\put(-15,25){\makebox(0,0){Donaldson-Thomas}}
\put(-15,15){\makebox(0,0){theory of $\P ^{1}\times \cnums ^{2}$}}
\end{picture}
\end{center}
\begin{quote}
\scriptsize{ The four isomorphic theories described by
$\mathsf{GW}^*(0\vertline 0,0)_{\lambda\mu\nu}$. The southern and
eastern theories have parameter $u$, while the northern and western
have parameter $q=-e^{iu}$. The horizontal equivalence is the
equivariant DT/GW correspondence for $\P ^{1} \times \cnums ^{2}$. The
vertical equivalence is the equivariant Crepant Resolution Conjecture
for $\Hilb \cnums ^{2}\to \operatorname{Sym}\cnums ^{2}$. The four
corners are computed in \cite{Bryan-Graber,dtlc,Ok-Pan-Hilb} and the
present paper.  }
\end{quote}

\normalsize

\pagebreak

\appendix 
\section{Appendix: Reconstruction Result}\label{appendix}



{By \large J. Bryan, C. Faber, A. Okounkov, and R. Pandharipande}



\subsection{Overview}

We present a proof of Theorem~\ref{thm: reconstruction of level (0,0)
pants} using a closed formula for the series
$$\mathsf{GW}^*(0\vertline 0,0)_{\lambda,(2),\nu}$$
obtained from Theorem \ref{thm: the dd2 integral}
and the semisimplicity of the Frobenius algebra associated to the
level $(0,0)$ theory. The proof was 
motivated by the study of the quantum cohomology of $\Hilb ^{n} (\cnums
^{2})$ in \cite{Ok-Pan-Hilb}.

\subsection{Fock space}


By definition, the Fock space $\cF$ 
is freely generated over $\Q$ by commuting 
creation operators $$\alpha_{-k},\ \ k\in\Z_{>0},$$
acting on the vacuum vector $\vac$. The annihilation 
operators 
$$\alpha_{k},\ \  k\in\Z_{>0},$$ kill the vacuum 
$$
\alpha_k \cdot \vac =0,\quad k>0 \,,
$$
and satisfy the commutation relations
$$
\left[\alpha_k,\alpha_l\right] = k \, \delta_{k+l,0}\,. 
$$

A natural basis of $\cF$ is given by 
the vectors  
\begin{equation}
  \label{basis}
  \lv \mu \rang = \frac{1}{\zz(\mu)} \, \prod_{i=1}^{l (\mu )} \alpha_{-\mu_i} \, \vac \,.
\end{equation}
indexed by partitions 
$\mu$. 
After extending scalars to $\Q(t_1,t_2)$,
we define the following 
{\em nonstandard} inner product on $\cF$:
\begin{equation}
  \label{inner_prod}
  \lang \mu | \nu \rang = 
\frac{(-1)^{|\mu|-\ell(\mu)}}{(t_1 t_2)^{l(\mu)}} 
\frac{\delta_{\mu\nu}}{\zz(\mu)} \,. 
\end{equation}

\subsection{The matrix $\MM _{2}$}

Define the linear transformation $\MM_2$ on $\cF$ by
\[
\lang \mu  | \MM_2 | \nu \rang = (-1)^{|\mu|}
\mathsf{GW}^*(0\vertline 0,0)_{\mu, (2), \nu}\  \delta_{|\mu|,|\nu|},
\]
after an extension of scalars to $\Q(t_1,t_2)[[u]]$.

The matrix $\MM_2$ can be written in closed form in terms of creation and
annihilation operators on Fock space:
\begin{multline}
\label{vv22}
-\MM_2 = \frac{t_1+t_2}{2} \sum_{k>0} \left( {k}  \frac{(-q)^k+1}{(-q)^k-1} -\frac{(-q)+1}{(-q)-1}\right) \, \alpha_{-k} \, \alpha_k  + \\
\frac12 \sum_{k,l>0} 
\Big[t_1 t_2 \, \alpha_{k+l} \, \alpha_{-k} \, \alpha_{-l} -
 \alpha_{-k-l}\,  \alpha_{k} \, \alpha_{l} \Big] \,,
\end{multline}
where $-q=e^{iu}$.  
The above formula was studied in \cite{Ok-Pan-Hilb}  as
the matrix of quantum multiplication by the hyperplane class
in the quantum cohomology of $\Hilb^n(\cnums^2)$.

Formula~\eqref{vv22} can be obtained as follows. Using dimension counts
similar to those in subsection~\ref{subsec: level 0 tube and cap}, the
disconnected invariants $\mathsf{GW}^* (0\vertline 0,0)_{\mu, (2),\nu }$ 
are easily reduced to connected invariants of one of two
possible types. First, there are the (necessarily connected) invariants
$\mathsf{GW}^*(0\vertline 0,0)_{(d) ,(2), (d)}$ (computed in
Theorem~\ref{thm: the dd2 integral}), and second there are domain genus 0
Hurwitz numbers. The combinatorics of writing disconnected invariants in
terms of connected invariants is most efficiently handled with the Fock
space formalism and yields Equation~\eqref{vv22}.

The first summand of Equation~\eqref{vv22} gives the diagonal terms of the
matrix $\MM _{2}$.  The second summand gives the off diagonal terms with
the $t_{1}t_{2}$ term of the summand appearing below the diagonal and the
remaining term appearing above.

\begin{lemma}\label{lem: eigenvals of M2 are distinct}
The eigenvalues of $\MM _{2} $ are distinct.
\end{lemma}

The eigenvalues are symmetric functions in $t_{1}$ and $t_{2}$. In the
$t_{1}t_{2}=0$ limit, $\MM _{2}$ is \emph{upper-triangular}. Hence it
suffices to show that the diagonal entries are distinct. By
Equation~\eqref{vv22}, the diagonal entry at a partition $\mu $ is
\begin{equation}\label{eqn: diagonal entries of M}
-\frac{t_{1}+t_{2}}{2}\sum _{k>0}km_{k} (\mu )F_{k}
\end{equation}
where $m_{k} (\mu )$ is the number of $k$'s in the partition $\mu $ and
\[
F_{k}=k\frac{(-q)^{k}+1}{(-q)^{k}-1}-\frac{(-q)+1}{(-q)-1}.
\]
The rational functions $\{F_{k} \}_{k>1}$ are easily seen to be linearly
independent over $\qnums $ (by, for example, studying the poles of
$F_{k}$), and hence  the diagonal entries  are distinct.\qed

\subsection{Proof of Theorem \ref{thm: reconstruction of level (0,0) pants}}

Let $d>0$.  We
abbreviate a list $(2),\dots ,(2)$ of $r$ copies of $(2)$ by $(2)^{r}$.
The gluing formula yields the equation
\[
\lang \mu  \vertline  \MM^r_2 \vertline  \nu \rang 
 =
(-1)^d\mathsf{GW}^*(0\vertline 0,0)_{\mu, (2)^{r}, \nu}
\]
for partitions $\mu,\nu$ of $d$.

A second application of the gluing formula yields the following computation:
\begin{eqnarray*}
\mathsf{GW}^*(0\vertline 0,0)_{\mu, (2)^r, \nu} & = & 
\sum_{\gamma\vdash d} \mathsf{GW}^*(0\vertline 0,0)_{\mu\gamma\nu} 
\mathsf{GW}^*(0\vertline 0,0)^{\gamma }_{ (2)^r} \\
& = & 
\sum_{\gamma\vdash d}
\mathsf{GW}^*(0\vertline 0,0)_{\mu\gamma\nu} 
\ \zz(\gamma)(-t_1t_2)^{l(\gamma)}  (-1)^d  \lang \gamma  \vertline \MM^r_2 \vertline (1^d) \rang.
\end{eqnarray*}
The second equality uses the level $(0,0)$ cap calculation of
Lemma~\ref{lem: charge (0,0) cap}.

Taken together, the above equation provide a linear system for the
degree $d$, level $(0,0)$ pair of pants integrals,
\begin{equation}\label{ff897}
\lang \mu  \vertline \MM^r_2 \vertline \nu \rang = 
\sum_{\gamma\vdash d}
\mathsf{GW}^*(0\vertline 0,0)_{\mu\gamma\nu} 
\ \zz(\gamma)(-t_1t_2)^{l(\gamma)}  
\  \lang \gamma  \vertline \MM^r_2 \vertline (1^d) \rang.
\end{equation}
The linear equations have coefficients in the field $\Q(t_1,t_2,q)$.
The proof of the Theorem is concluded by demonstrating the nonsingularity of
the system \eqref{ff897}.

Let $\cF_d \subset \cF$ be the subspace spanned by the vectors $|\mu
\rangle$ satisfying $|\mu|=d$. The transformation $\MM_2$ preserves
$\cF_d$.  

The eigenvectors for $\MM _{2}$ restricted to $\cF _{d}$ are the idempotent
basis of the semisimple Frobenius algebra associated to the degree $d$,
level $(0,0)$ theory.  The identity vector $|(1^d)\rangle$ of the Frobenius
algebra has the coefficient 1 in each component of the idempotent
basis. Hence, the set of vectors
\[
\{ \ \MM_2^r\vertline (1^d) \rangle \ \}_{r\geq 0}
\]
have coefficients given by powers of the eigenvalues of $\MM _{2}$
restricted to $\cF _{d}$. These eigenvalues are distinct by Lemma~\ref{lem:
eigenvals of M2 are distinct}, thus the above set of vectors spans $\cF_d$
and the linear system \eqref{ff897} is nonsingular. \qed

\subsection{Proof of Theorem \ref{thm: Z is a rational function of q}}

Since $-q=e^{-iu}$, we may disregard all integral terms in the exponent of the
prefactor 
$$e^{\frac{idu}{2} (2-2g+{k_{1}}+{k_{2}})}.$$
Consider
the product
\[
e^{\frac{idu}{2} ({k_{1}}+{k_{2}})}\  
\mathsf{GW}^*(g\vertline  k_{1},k_{2})_{\lambda ^{1}\dots \lambda ^{r}}.
\]
The series $\mathsf{GW}^*(g\vertline k_{1},k_{2})_{\lambda ^{1}\dots
\lambda ^{r}}$ can be calculated by gluing in terms of the caps
\begin{equation}\label{xsd} 
\mathsf{GW}^*
(0\vertline \pm1,0)_\lambda, \ \ \mathsf{GW}^*(0\vertline 0,\pm1)_{\lambda }
\end{equation}
and the pair of pant series
\[
\mathsf{GW}^*(0\vertline 0,0)_{\lambda\mu\nu}.
\]
By the proof of Theorem 
\ref{thm: reconstruction of level (0,0) pants}, 
the pair of pant series lie in $\Q(t_1,t_2,q)$.
By the calculation of Section \ref{n3451},
$$e^{-\frac{idu}{2}}\mathsf{GW}^*(0\vertline -1,0)_\lambda, \
e^{-\frac{idu}{2}}\mathsf{GW}^*(0\vertline 0,-1)_\lambda \in
\Q(t_1,t_2,q).$$ Since the opposite caps are inverses in the Frobenius
algebra, we conclude
$$e^{\frac{idu}{2}}\mathsf{GW}^*(0\vertline 1,0)_\lambda, \
e^{\frac{idu}{2}}\mathsf{GW}^*(0\vertline 0,1)_\lambda \in \Q(t_1,t_2,q).$$
The Theorem is proven by distributing a factor of $e^{\pm\frac{idu}{2}}$ to
each cap of type \eqref{xsd} in the gluing formula. \qed


\bibliographystyle{plain}

\vspace{+10 pt}
\noindent
Department of Mathematics \\
University of British Columbia \\
Vancouver, BC, V6T 1Z2, Canada\\
jbryan@math.ubc.ca \\

\vspace{+10 pt}
\noindent
Department of Mathematics\\
Princeton University\\
Princeton, NJ 08544, USA\\
rahulp@math.princeton.edu

\end{document}